\documentclass[12pt]{article}
\usepackage{amssymb,amsmath,latexsym,amsthm}
\usepackage{url}
\usepackage{enumerate}
\usepackage{ifpdf}
\ifpdf
  \usepackage[pdftex]{hyperref}
\else
  \usepackage[dvips]{hyperref}
\fi

\usepackage{endnotes}

\title{The weak Stratonovich integral with respect to fractional Brownian motion with Hurst parameter $1/6$}
\author{
  Ivan Nourdin\thanks{Supported in part by the (french) ANR grant `Exploration des Chemins Rugueux'.}\\
  Universit\'e Paris VI
  \and
  Anthony R\'eveillac\thanks{Supported by DFG research center Matheon project E2.}\\
  Humboldt-Universit\"at zu Berlin
  \and
  Jason Swanson\thanks{Supported in part by NSA grant H98230-09-1-0079.}\\
  University of Central Florida
  }

\DeclareMathOperator{\sgn}{sgn}

\DeclareMathOperator{\tot}{\wt\ot}
\DeclareMathOperator{\dom}{Dom}
\usepackage{color,eufrak,mathrsfs}

\begin{document}

\newtheorem{thm}{Theorem}[section]
\newtheorem{corollary}[thm]{Corollary}
\newtheorem{prop}[thm]{Proposition}
\newtheorem{lemma}[thm]{Lemma}
\theoremstyle{definition}
\newtheorem{remark}[thm]{Remark}

\numberwithin{equation}{section}

\def\al{\alpha}
\def\be{\beta}
\def\ga{\gamma}
\def\de{\delta}
\def\De{\Delta}
\def\ep{\varepsilon}
\def\th{\theta}
\def\ka{\kappa}
\def\la{\lambda}
\def\si{\sigma}
\def\ph{\varphi}
\def\om{\omega}
\def\Om{\Omega}

\def\AA{\mathcal{A}}
\def\AAA{\mathbb{A}}
\def\CC{\mathbb{C}}
\def\FF{\mathcal{F}}
\def\GG{\mathcal{G}}
\def\II{\mathcal{I}}
\def\NN{\mathbb{N}}
\def\PP{\mathcal{P}}
\def\RR{\mathbb{R}}
\def\VV{\mathbb{V}}
\def\WW{\mathbb{W}}
\def\ZZ{\mathbb{Z}}

\def\wt{\widetilde}
\def\wh{\widehat}
\def\ot{\otimes}

\def\HH{\EuFrak H}
\def\sk{{\mathbb{D}}}
\def\lcr{\left[}
\def\rcr{\right]}
\newcommand{\R}{\mathbb{R}}
\def\real{\mathbb{R}}

\def\pa{\partial}
\def\ds{\displaystyle}
\def\ts{\textstyle}

\def\pf{\noindent{\bf Proof.} }
\def\qed{\hfill $\Box$}

\providecommand{\flr}[1]{\lfloor#1\rfloor}
\providecommand{\ang}[1]{\langle#1\rangle}
\providecommand{\cub}[1]{[\![#1]\!]}

\maketitle

\begin{abstract}

Let $B$ be a fractional Brownian motion with Hurst parameter $H=1/6$. It is known that the symmetric Stratonovich-style Riemann sums for $\int g(B(s))\,dB(s)$ do not, in general, converge in probability. We show, however, that they do converge in law in the Skorohod space of c\`adl\`ag functions. Moreover, we show that the resulting stochastic integral satisfies a change of variable formula with a correction term that is an ordinary It\^o integral with respect to a Brownian motion that is independent of $B$.

\bigskip


\bigskip


\end{abstract}

\section{Introduction}

The Stratonovich integral of $X$ with respect to $Y$, denoted $\int_0^tX(s) \circ dY(s)$, can be defined as the limit in probability, if it exists, of
  \begin{equation}\label{strat_sum}
  \sum_{t_j\le t}\frac{X(t_{j-1}) + X(t_j)}2(Y(t_j) - Y(t_{j-1})),
  \end{equation}
as the mesh of the partition $\{t_j\}$ goes to zero. Typically, we regard \eqref{strat_sum} as a process in $t$, and require that it converges uniformly on compacts in probability (ucp).

This is closely related to the so-called symmetric integral, denoted by $\int_0^t X(s)\,d^\circ Y(s)$, which is the ucp limit, if it exists, of
  \begin{equation}\label{symm_int}
  \frac1\ep\int_0^t\frac{X(s) + X(s+\ep)}2(Y(s+\ep) - Y(s))\,ds,
  \end{equation}
as $\ep\to 0$. The symmetric integral is an example of the regularization procedure, introduced by Russo and Vallois, and on which there is a wide body of literature. For further details on stochastic calculus via regularization, see the excellent survey article \cite{RV} and the many references therein.

A special case of interest that has received considerable attention in the literature is when $Y=B^H$, a fractional Brownian motion with Hurst parameter $H$. It has been shown independently in \cite{CN} and \cite{GNRV} that when $Y=B^H$ and $X=g(B^H)$ for a sufficiently differentiable function $g(x)$, the symmetric integral exists for all $H>1/6$. Moreover, in this case, the symmetric integral satisfies the classical Stratonovich change of variable formula,
  \[
  g(B^H(t)) = g(B^H(0)) + \int_0^t g'(B^H(s))\,d^\circ B^H(s).
  \]
However, when $H=1/6$, the symmetric integral does not, in general, exist. Specifically, in \cite{CN} and \cite{GNRV}, it is shown that \eqref{symm_int} does not converge in probability when $Y=B^{1/6}$ and $X=(B^{1/6})^2$. It can be similarly shown that, in this case, \eqref{strat_sum} also fails to converge in probability.

This brings us naturally to the notion which is the focus of this paper: the weak Stratonovich integral, which is the limit in law, if it exists, of \eqref{strat_sum}. We focus exclusively on the case $Y=B^{1/6}$. For simplicity, we omit the superscript and write $B=B^{1/6}$. Our integrands shall take the form $g(B(t))$, for $g\in C^\infty(\RR)$, and we shall work only with the uniformly spaced partition, $t_j=j/n$. In this case, \eqref{strat_sum} becomes
  \[
  I_n(g,B,t) = \sum_{j=1}^{\flr{nt}}
    \frac{g(B(t_{j-1})) + g(B(t_j))}2\De B_j,
  \]
where $\flr{x}$ denotes the greatest integer less than or equal to $x$, and $\De B_j=B(t_j)-B(t_{j-1})$. We show that the processes $I_n(g,B)$ converge in law in $D_\RR[0,\infty)$, the Skorohod space of c\`adl\`ag functions from $[0,\infty)$ to $\RR$. We let $\int_0^t g(B(s))\,dB(s)$ denote a process with this limiting law, and refer to this as the weak Stratonovich integral.

The weak Stratonovich integral with respect to $B$ does not satisfy the classical Stratonovich change of variable formula. Rather, we show that it satisfies a change of variable formula with a correction term that is a classical It\^o integral. Namely,
  \begin{equation}\label{ito_form}
  g(B(t)) = g(B(0)) + \int_0^t g'(B(s))\,dB(s)
    - \frac1{12}\int_0^t g'''(B(s))\,d\cub{B}_s,
  \end{equation}
where $\cub{B}$ is what we call the signed cubic variation of $B$. That is, $\cub{B}$ is the limit in law of the sequence of processes $V_n(B,t)=\sum_{j=1}^{\flr{nt}}\De B_j^3$. It is shown in \cite{NO} that $\cub{B}=\ka W$, where $W$ is a standard Brownian motion, independent of $B$, and $\ka\simeq 2.322$. (See \eqref{kappa} for the exact definition of $\ka$.) The correction term in \eqref{ito_form} is then a standard It\^o integral with respect to Brownian motion.

Our precise results are actually somewhat stronger than this, in that we prove the joint convergence of the processes $B$, $V_n(B)$, and $I_n(g,B)$. (See Theorem \ref{T:main}.) We also discuss the joint convergence of multiple sequences of Riemann sums for different integrands. (See Theorem \ref{T:joint_conv} and Remark \ref{R:joint_conv}.)

The work in this paper is a natural follow-up to \cite{BS} and \cite{NR}. There, analogous results were proven for $B^{1/4}$ in the context of midpoint-style Riemann sums. The results in \cite{BS} and \cite{NR} were proven through different methods, and in the present work, we combine the two approaches to prove our main results.

Finally, let us stress the fact that, as a byproduct of the proof of \eqref{ito_form}, we show in the present paper that
  \[
  n^{-1/2}\sum_{j=1}^{\flr{n\cdot}}g(B(t_{j-1}))
    \,h_3(n^{1/6}\De B_j)
    \to -\frac18\int_0^\cdot g'''(B(s))\,ds
    + \int_0^\cdot g(B(s))\,d\cub{B}_s,
  \]
in the sense of finite-dimensional distributions on $[0,\infty)$, where $h_3(x) =x^3-3x$ denotes the third Hermite polynomial. (See more precisely Theorem \ref{fdd2} below. Also see Theorem \ref{T:fdd3}.) From our point of view, this result has also its own interest, and should be compared with the recent results obtained in \cite{NN,NNT}, concerning the weighted Hermite variations of fractional Brownian motion.

\section{Notation, preliminaries, and main result}

Let $B=B^{1/6}$ be a fractional Brownian motion with Hurst parameter $H=1/6$. That is, $B$ is a centered Gaussian process, indexed by $t\ge0$, such that
  \[
  R(s,t) = E[B(s)B(t)] = \frac12(t^{1/3} + s^{1/3} - |t - s|^{1/3}).
  \]
Note that $E|B(t)-B(s)|^2=|t-s|^{1/3}$. For compactness of notation, we will sometimes write $B_t$ instead of $B(t)$. Given a positive integer $n$, let $\De t=n^{-1}$ and $t_j=t_{j,n}=j\De t$. We shall frequently have occasion to deal with the quantity $\be_{j,n}=\be_j=(B(t_{j-1})+B(t_j))/2$. In estimating this and similar quantities, we shall adopt the notation $r_+=r\vee1$, which is typically applied to nonnegative integers $r$. We shall also make use of the Hermite polynomials,
  \begin{equation}\label{hermite}
  h_n(x) = (-1)^n e^{x^2/2}\frac{d^n}{dx^n}(e^{-x^2/2}).
  \end{equation}
Note that the first few Hermite polynomials are $h_0(x)=1$, $h_1(x)=x$, $h_2(x) =x^2-1$, and $h_3(x)=x^3-3x$. The following orthogonality property is well-known: if $U$ and $V$ are jointly normal with $E(U)=E(V)=0$ and $E(U^2)=E(V^2)=1$, then
  \begin{equation}\label{eq:ortho}
  E[h_p(U)h_q(V)] = \begin{cases}
    q!(E[UV])^q &\text{if $p=q$},\\
    0 &\text{otherwise}.
    \end{cases}
  \end{equation}

If $X$ is a c\`adl\`ag process, we write $X(t-)=\lim_{s\uparrow t}X(s)$ and $\De X(t)=X(t)-X(t-)$. The step function approximation to $X$ will be denoted by $X_n(t)=X(\flr{nt}/n)$, where $\flr{\cdot}$ is the greatest integer function. In this case, $\De X_n(t_{j,n})=X(t_j)-X(t_{j-1})$. We shall frequently use the shorthand notation $\De X_j=\De X_{j,n}=\De X_n(t_{j,n})$. For simplicity, positive integer powers of $\De X_j$ shall be written without parentheses, so that $\De X_j^k=(\De X_j)^k$.

The discrete $p$-th variation of $X$ is defined as
  \[
  V_n^p(X,t) = \sum_{j=1}^{\flr{nt}}|\De X_j|^p,
  \]
and the discrete signed $p$-th variation of $X$ is
  \[
  V_n^{p\pm}(X,t) = \sum_{j=1}^{\flr{nt}}|\De X_j|^p\sgn(\De X_j).
  \]
For the discrete signed cubic variation, we shall omit the superscript, so that
  \begin{equation}\label{vn3}
  V_n(X,t) = V_n^{3\pm}(X,t) = \sum_{j=1}^{\flr{nt}}\De X_j^3.
  \end{equation}
When we omit the index $t$, we mean to refer to the entire process. So, for example, $V_n(X)=V_n(X,\cdot)$ refers to the c\`adl\`ag process which maps $t\mapsto V_n(X,t)$.

Let $\{\rho(r)\}_{r\in\ZZ}$ be the sequence defined by
  \begin{equation}\label{rho}
  \rho(r) = \frac12(|r + 1|^{1/3} + |r - 1|^{1/3} - 2|r|^{1/3}).
  \end{equation}
Note that $\sum_{r\in\ZZ}|\rho(r)|<\infty$ and $E[\De B_i\De B_j] =n^{-1/3}\rho(i-j)$ for all $i,j\in\NN$. Let $\ka>0$ be defined by
  \begin{equation}\label{kappa}
  \ka^2 = 6\sum_{r\in\ZZ}\rho^3(r)
    = \frac34\sum_{r\in\ZZ}
    (|r + 1|^{1/3} + |r - 1|^{1/3} - 2|r|^{1/3})^3
    \simeq 5.391,
  \end{equation}
and let $W$ be a standard Brownian motion, defined on the same probability space as $B$, and independent of $B$. Define $\cub{B}_t=\ka W(t)$. We shall refer to the process $\cub{B}$ as the signed cubic variation of $B$. The use of this term is justified by Theorem \ref{T:signed_cubic}.

A function $g:\RR^d\to\RR$ has {\it polynomial growth} if there exist positive constants $K$ and $r$ such that $|g(x)| \le K(1 + |x|^r)$ for all $x\in\RR^d$. If $k$ is a nonnegative integer, we shall say that a function $g$ has {\it polynomial growth of order $k$} if $g\in C^k(\RR^d)$ and there exist positive constants $K$ and $r$ such that $|\pa^\al g(x)| \le K(1 + |x|^r)$ for all $x\in\RR^d$ and all $|\al|\le k$. (Here, $\al\in\NN_0^d=(\NN\cup\{0\})^d$ is a multi-index, and we adopt the standard multi-index notation: $\pa_j=\pa/\pa x_j$, $\pa^\al=\pa_1^{\al_1}\cdots\pa_d^{\al_d}$, and $|\al|=\al_1+\cdots+\al_d$.)

Given $g:\RR\to\RR$ and a stochastic process $\{X(t):t\ge0\}$, the Stratonovich Riemann sum will be denoted by
  \[
  I_n(g,X,t) = \sum_{j=1}^{\flr{nt}}
    \frac{g(X(t_{j-1})) + g(X(t_j))}2\De X_j.
  \]
The phrase ``uniformly on compacts in probability" will be abbreviated ``ucp." If $X_n$ and $Y_n$ are c\`adl\`ag processes, we shall write $X_n\approx Y_n$ or $X_n(t)\approx Y_n(t)$ to mean that $X_n-Y_n\to0$ ucp. In the proofs in this paper, $C$ shall denote a positive, finite constant that may change value from line to line.

\subsection{Conditions for relative compactness}

The Skorohod space of c\`adl\`ag functions from $[0,\infty)$ to $\RR^d$ is denoted by $D_{\RR^d}[0,\infty)$. Note that $D_{\RR^d}[0,\infty)$ and $(D_\RR[0,\infty))^d$ are not the same. In particular, the map $(x,y)\mapsto x+y$  is continuous from $D_{\RR^2}[0,\infty)$ to $D_\RR[0,\infty)$, but it is not continuous from $(D_\RR[0,\infty))^2$ to $D_\RR[0,\infty)$. Convergence in $D_{\RR^d} [0,\infty)$ implies convergence in $(D_\RR[0,\infty))^d$, but the converse is not true.

Note that if the sequences $\{X_n^{(1)}\},\ldots,\{X_n^{(d)}\}$ are all relatively compact in $D_\RR[0,\infty)$, then the sequence of $d$-tuples $\{(X_n^{(1)},\ldots, X_n^{(d)})\}$ is relatively compact in $(D_\RR[0,\infty))^d$. It may not, however, be relatively compact in $D_{\RR^d}[0,\infty)$. We will therefore need the following well-known result. (For more details, see Section 2.1 of \cite{BS} and the references therein.)

\begin{lemma}\label{L:prod_space}
Suppose $\{(X_n^{(1)},\ldots,X_n^{(d)})\}_{n=1}^\infty$ is relatively compact in $(D_\RR [0,\infty))^d$. If, for each $j\ge 2$, the sequence $\{X_n^{(j)}\}_ {n=1} ^\infty$ converges in law in $D_\RR[0,\infty)$ to a continuous process, then $\{(X_n^{(1)},\ldots, X_n^{(d)})\}_{n=1}^\infty$ is relatively compact in $D_{\RR^d}[0, \infty)$.
\end{lemma}

Our primary criterion for relative compactness is the following moment condition, which is a special case of Corollary 2.2 in \cite{BS}.

\begin{thm}\label{T:moment_cond}
Let $\{X_n\}$ be a sequence of processes in $D_{\RR^d}[0,\infty)$. Let $q(x) = |x|\wedge1$. Suppose that for each $T>0$, there exists $\nu>0$, $\be>0$, $C>0$, and $\th>1$ such that $\sup_n E[|X_n(T)|^\nu]<\infty$ and
  \begin{equation}\label{moment_cond}
  E[q(X_n(t) - X_n(s))^\be]
    \le C\bigg(\frac{\flr{nt} - \flr{ns}}n\bigg)^\th,
  \end{equation}
for all $n$ and all
$0\le s\le t\le T$.
 Then $\{X_n\}$ is relatively compact.
\end{thm}

Of course, a sequence $\{X_n\}$ converges in law in $D_{\RR^d}[0,\infty)$ to a process $X$ if $\{X_n\}$ is relatively compact and $X_n\to X$ in the sense of finite-dimensional distributions on $[0,\infty)$. We shall also need the analogous theorem for convergence in probability, which is Lemma A2.1 in \cite{DK}. Note that if $x:[0,\infty)\to\RR^d$ is continuous, then $x_n\to x$ in $D_{\RR^d}[0,\infty)$ if and only if $x_n\to x$ uniformly on compacts.

\begin{lemma}\label{L:in_prob}
Let $\{X_n\},X$ be processes with sample paths in $D_{\RR^d}[0,\infty)$ defined on the same probability space. Suppose that $\{X_n\}$ is relatively compact in $D_{\RR^d} [0,\infty)$ and that for a dense set $H\subset[0,\infty)$, $X_n(t)\to X(t)$ in probability for all $t\in H$. Then $X_n\to X$ in probability in $D_{\RR^d} [0,\infty)$. In particular, if $X$ is continuous, then $X_n\to X$ ucp.
\end{lemma}

We will also need the following lemma, which is easily proved using the Prohorov metric.

\begin{lemma}\label{L:truncate}
Let $(E,r)$ be a complete and separable metric space. Let $X_n$ be a sequence of $E$-valued random variables and suppose, for each $k$, there exists a sequence $\{X_{n,k}\}_{n=1}^\infty$ such that $\limsup_{n\to\infty}E[r(X_n, X_{n,k})] \le \de_k$, where $\de_k\to0$ as $k\to\infty$. Suppose also that for each $k$,
there exists $Y_k$ such that $X_{n,k}\to Y_k$ in law as $n\to\infty$.
Then there exists $X$ such that $X_n\to X$ in law and $Y_k\to X$ in law.
\end{lemma}

\subsection{Elements of Malliavin calculus}\label{prelim}

In the sequel, we will need some elements of Malliavin calculus that we collect here. The reader is referred to \cite{Malliavin} or \cite{Nualart} for any unexplained notion discussed in this section.

We denote by $X=\{X(\ph):\,\ph\in\HH\}$ an isonormal Gaussian process over $\HH$, a real and separable Hilbert space. By definition, $X$ is a centered Gaussian family indexed by the elements of $\HH$ and such that, for every $\ph,\psi\in\HH$,
  \[
  E[X(\ph)X(\psi)]=\ang{\ph,\psi}_\HH.
  \]
We denote by $\HH^{\ot q}$ and $\HH^{\odot q}$, respectively, the tensor space and the symmetric tensor space of order $q\ge 1$. Let $\mathscr{S}$ be the set of cylindrical functionals $F$ of the form
  \begin{equation}\label{eq:cylindrical}
  F = f(X(\ph_1),\ldots,X(\ph_n)),
  \end{equation}
where $n\ge 1$, $\ph_i \in \HH$ and the function $f\in C^{\infty} (\RR^n)$ is such that its partial derivatives have polynomial growth. The Malliavin derivative $DF$ of a functional $F$ of the form \eqref{eq:cylindrical} is the square integrable $\HH$-valued random variable defined as
  \[
  DF = \sum_{i=1}^n
    \frac{\pa f}{\pa x_i}(X(\ph_1),\ldots,X(\ph_n))\ph_i.
  \]
In particular, $DX(\ph) = \ph$ for every $\ph\in \HH$. By iteration, one can define the $m$th derivative $D^m F$ (which is an element of $L^2(\Omega, \HH^{\odot m})$) for every $m\ge 2$, giving
  \[
  D^mF = \sum_{i_1,\ldots,i_m}^n \frac{\pa^m f}
    {\pa x_{i_1}\cdots\pa x_{i_m}}
    (X(\ph_1),\ldots,X(\ph_n))
    \ph_{i_1}\otimes\cdots\otimes\ph_{i_m}.
  \]
As usual, for $m\ge 1$, $\sk^{m,2}$ denotes the closure of $\mathscr{S}$ with respect to the norm $\|\cdot\|_{m,2}$, defined by the relation
  \[
  \|F\|_{m,2}^2 = EF^2 + \sum_{i=1}^m
    E\|D^i F\|_{\HH^{\ot i}}^2.
  \]
The Malliavin derivative $D$ satisfies the following chain rule: if $f:\R^n\rightarrow\R$ is in $C^1_b$ (that is, the collection of continuously differentiable functions with a bounded derivative) and if $\{F_i\}_{i=1,\ldots,n}$ is a vector of elements of $\sk^{1,2}$, then $f(F_1,\ldots,F_n)\in\sk^{1,2}$ and
  \begin{equation}\label{chainrule}
  Df(F_1,\ldots,F_n) = \sum_{i=1}^n
    \frac{\pa f}{\pa x_i}(F_1,\ldots, F_n)DF_i.
  \end{equation}
This formula can be extended to higher order derivatives as
  \begin{equation}\label{m-chainrule}
  D^mf(F_1,\ldots,F_n) = \sum_{v\in\PP_m}C_v
    \sum_{i_1,\ldots,i_k=1}^n
    \frac{\pa^k f}{\pa x_{i_1}\cdots\pa x_{i_k}}
    (F_1,\ldots,F_n)D^{v_1}F_{i_1}\tot
    \cdots\tot D^{v_k}F_{i_k},
  \end{equation}
where $\PP_m$ is the set of vectors $v=(v_1,\ldots,v_k)\in\NN^k$ such that $k\ge1$, $v_1\le\cdots\le v_k$, and $v_1+\cdots+v_k=m$. The constants $C_v$ can be written explicitly as $C_v=m!(\prod_{j=1}^n m_j!(j!)^{m_j})^{-1}$, where $m_j=|\{\ell:v_\ell=j\}|$.

\begin{remark}\label{rem-lei-rule}
In \eqref{m-chainrule}, $a\tot b$ denotes the symmetrization of the tensor product $a\ot b$. Recall that, in general, the {\it symmetrization} of a function $f$ of $m$ variables is the function $\wt f$ defined by
  \begin{equation}\label{symmetrization}
  \wt f(t_1,\ldots,t_m) = \frac1{m!}
    \sum_{\si\in\mathfrak{S}_m}f(t_{\si(1)},\ldots,t_{\si(m)}),
  \end{equation}
where $\mathfrak{S}_m$ denotes the set of all permutations of $\{1, \ldots,m\}$.
\end{remark}

We denote by $I$ the adjoint of the operator $D$, also called the divergence operator. A random element $u \in L^2(\Omega, \HH)$ belongs to the domain of $I$, noted $\dom(I)$, if and only if it satisfies
  \[
  |E\ang{DF,u}_\HH|\le c_u\sqrt{EF^2}\quad
    \text{for any $F\in\mathscr{S}$},
  \]
where $c_u$ is a constant depending only on $u$. If $u \in \dom(I)$, then the random variable $I(u)$ is defined by the duality relationship (customarily called ``integration by parts formula"):
  \begin{equation}\label{ipp}
  E[F I(u)] = E\ang{DF,u}_\HH,
  \end{equation}
which holds for every $F \in \sk^{1,2}$.

For every $n\geq 1$, let $\mathcal{H}_n$ be the $n$th Wiener chaos of $X$, that is, the closed linear subspace of $L^2$ generated by the random variables $\{h_n(X(\ph)):\ph\in \HH,|\ph|_\HH=1\}$, where $h_n$ is the Hermite polynomial defined by \eqref{hermite}. The mapping
  \begin{equation}\label{eq:defIn}
  I_n(\ph^{\ot n}) = h_n(X(\ph))
  \end{equation}
provides a linear isometry between the symmetric tensor product $\HH^{\odot n}$ (equipped with the modified norm $\frac1{\sqrt{n!}}\|\cdot\|_{\HH^{\ot n}}$) and $\mathcal{H}_n$. The following duality formula holds:
  \begin{equation}\label{dual}
  E[FI_n(f)] = E\ang{D^{n}F,f}_{\HH^{\ot n}},
  \end{equation}
for any element $f\in \HH^{\odot n}$ and any random variable $F\in \mathbb{D}^{n,2}$. We will also need the following particular case of the classical product formula between multiple integrals: if $\ph, \psi\in \HH$ and $m,n\ge 1$, then
  \begin{equation}\label{multiplication}
  I_m(\ph^{\ot m})I_n(\psi^{\ot n})
    = \sum_{r=0}^{m\wedge n} r!\binom mr\binom nr
    I_{m+n-2r}(\ph^{\ot(m-r)} \ot \psi^{\ot(n-r)})
    \ang{\ph,\psi}_\HH^r.
  \end{equation}

Finally, we mention that the Gaussian space generated by $B=B^{1/6}$ can be identified with an isonormal Gaussian process of the type $B=\{B(h):h\in\HH\}$, where the real and separable Hilbert space $\HH$ is defined as follows: (i) denote by $\mathscr{E}$ the set of all $\RR$-valued step functions on $[0,\infty)$, (ii) define $\HH$ as the Hilbert space obtained by closing $\mathscr{E}$ with respect to the scalar product
  \[
  \ang{{\bf 1}_{[0,t]},{\bf 1}_{[0,s]}}_\HH
    = E[B(s)B(t)] = \frac12(t^{1/3}+s^{1/3}-|t-s|^{1/3}).
  \]
In particular, note that $B(t)=B(\mathbf{1}_{[0,t]})$. To end up, let us stress that the $m$th derivative $D^m$ (with respect to $B$) verifies the Leibniz rule. That is, for any $F,G\in\mathbb{D}^{m,2}$ such that $FG\in\mathbb{D}^{m,2}$, we have
  \begin{equation}\label{leibnitz}
  D^m_{t_1,\ldots,t_m}(FG) = \sum D_J^{|J|}(F)D_{J^c}^{m-|J|}(G),
    \quad t_i\in[0,T],\quad i=1,\ldots,m,
  \end{equation}
where the sum runs over all subsets $J$ of $\{t_1,\ldots,t_m\}$, with $|J|$ denoting the cardinality of $J$. Note that we may also write this as
  \begin{equation}\label{tensor-leibnitz}
  D^m(FG) = \sum_{k=0}^m \binom mk
    (D^k F)\tot(D^{m-k}G).
  \end{equation}

\subsection{Expansions and Gaussian estimates}

A key tool of ours will be the following version of Taylor's theorem with remainder.

\begin{thm}\label{T:Taylor}
Let $k$ be a nonnegative integer. If $g\in C^k(\RR^d)$, then
  \[
  g(b) = \sum_{|\al|\le k} \pa^\al g(a)\frac{(b - a)^\al}{\al!}
    + R_k(a,b),
  \]
where
  \[
  R_k(a,b) = k\sum_{|\al|=k}\frac{(b - a)^\al}{\al!}
    \int_0^1 (1 - u)^k [\pa^\al g(a + u(b - a)) - \pa^\al g(a)]\,du
  \]
if $k\ge 1$, and $R_0(a,b)=g(b)-g(a)$. In particular, $R_k(a,b) = \sum_{|\al|=k} h_\al(a,b)(b - a)^\al$, where $h_\al$ is a continuous function with $h_\al(a,a)=0$ for all $a$. Moreover,
  \[
  |R_k(a,b)| \le (k \vee 1)\sum_{|\al|=k}M_\al|(b - a)^\al|,
  \]
where $M_\al=\sup\{|\pa^\al g(a + u(b - a))-\pa^\al g(a)|:0\le u\le 1\}$.
\end{thm}

The following related expansion theorem is a slight modification of Corollary 4.2 in \cite{BS}.

\begin{thm}\label{T:Gauss_Taylor}
Recall the Hermite polynomials $h_n(x)$ from \eqref{hermite}. Let $k$ be a nonnegative integer. Suppose $\ph:\RR\to\RR$ is measurable and has polynomial growth with constants $\wt K$ and $r$. Suppose $f\in C^{k+1}(\RR^d)$ has polynomial growth of order $k+1$, with constants $K$ and $r$. Let $\xi\in\RR^d$ and $Y\in\RR$ be jointly normal with mean zero. Suppose that $EY^2=1$ and $E\xi_j^2\le\nu$ for some $\nu>0$. Define $\eta\in\RR^d$ by $\eta_j=E[\xi_jY]$. Then
  \[
  E[f(\xi)\ph(Y)] = \sum_{|\al|\le k}\frac1{\al!}\,\eta^\al
    E[\pa^\al f(\xi)]E[h_{|\al|}(Y)\ph(Y)] + R,
  \]
where $|R|\le CK|\eta|^{k+1}$ and $C$ depends only on $\wt K$, $r$, $\nu$, $k$, and $d$.
\end{thm}

\pf Although this theorem is very similar to Corollary 4.2 in \cite{BS}, we provide here another proof by means of Malliavin calculus.

Observe first that, without loss of generality, we can assume that $\xi_i=X(v_i)$, $i=1,\ldots,d$, and $Y=X(v_{d+1})$, where $X$ is an isonormal process over $\HH=\R^{d+1}$ and where $v_1,\ldots,v_{d+1}$ are some adequate vectors belonging in $\HH$. Since $\ph$ has polynomial growth, we can expand it in terms of Hermite polynomials, that is $\ph=\sum_{q=0}^\infty c_q h_q$. Thanks to \eqref{eq:ortho}, note that $q!c_q=E[\ph(Y)h_q(Y)]$. We set
  \[
  \wh\ph_k = \sum_{q=0}^k c_q h_q\quad\text{and}\quad
    \check{\ph}_k = \sum_{q=k+1}^\infty c_q h_q.
  \]
Of course, we have
  \[
  E[f(\xi)\ph(Y)] = E[f(\xi)\wh\ph_k(Y)] + E[f(\xi)\check{\ph}_k(Y)].
  \]
We obtain
  \begin{align*}
  E[f(\xi)\wh\ph_k(Y)] &= \sum_{q=0}^k \frac1{q!}
    \,E[\ph(Y)h_q(Y)]\,E[f(\xi)h_q(Y)]\\
  &= \sum_{q=0}^k \frac1{q!}\,E[\ph(Y)h_q(Y)]
    \,E[f(\xi)I_q(v_{d+1}^{\ot q})]\quad\text{by \eqref{eq:defIn}}\\
  &= \sum_{q=0}^k \frac1{q!}\,E[\ph(Y)h_q(Y)]
    \,E[\ang{D^qf(\xi),v_{d+1}^{\ot q}}_{\HH^{\ot q}}]
    \quad\text{by \eqref{dual}}\\
  &= \sum_{q=0}^k \frac1{q!}\sum_{i_1,\ldots,i_q=1}^d
    \,E[\ph(Y)h_q(Y)]
    \,E\left[{\frac{\pa^qf}{\pa x_{i_1}\cdots\pa x_{i_q}}(\xi)}\right]
    \prod_{\ell=1}^q\eta_{i_\ell}
    \quad\text{by \eqref{m-chainrule}}.
  \end{align*}
Since the map $\Phi:\{1,\ldots,d\}^q\to\{\al\in\NN_0^d:|\al|=q\}$ defined by $(\Phi(i_1,\ldots,i_q))_j=|\{\ell:i_\ell=j\}|$ is a surjection with $|\Phi^{-1}(\al)|=q!/\al!$, this gives
  \begin{align*}
  E[f(\xi)\wh\ph_k(Y)] &= \sum_{q=0}^k \frac1{q!}\sum_{|\al|=q}
    \frac{q!}{\al!}\,E[\ph(Y)h_q(Y)]
    \,E[\pa^\al f(\xi)]\eta^\al\\
  &= \sum_{|\al|\le k}\frac1{\al!}\,E[\ph(Y)h_{|\al|}(Y)]
    \,E[\pa^\al f(\xi)]\eta^\al.
  \end{align*}
On the other hand, the identity \eqref{eq:ortho}, combined with the fact that each monomial $x^n$ can be expanded in terms of the first $n$ Hermite polynomials, implies that $E[Y^{|\al|}\check{\ph}_k(Y)]=0$ for all $|\al|\le k$. Now, let $U=\xi-\eta Y$ and define $g:\RR^d\to\R$ by $g(x) =E[f(U+xY)\check{\ph}_k(Y)]$. Since $\ph$ (and, consequently, also $\check{\ph}_k$) and $f$ have polynomial growth, and all derivatives of $f$ up to order $k+1$ have polynomial growth, we may differentiate under the expectation and conclude that $g\in C^{k+1}(\RR^d)$. Hence, by Taylor's theorem (more specifically, by the version of Taylor's theorem which appears as Theorem 2.13 in \cite{BS}), and the fact that $U$ and $Y$ are independent,
  \begin{align*}
  E[f(\xi)\check{\ph}_k(Y)] = g(\eta) &= \sum_{|\al|\le k}
    \frac1{\al!}\eta^{\al}\pa^\al g(0) + R\\
  &= \sum_{|\al|\le k}\frac1{\al!}\eta^\al
    E[\pa^\al f(U)]E[Y^{|\al|}\check{\ph}_k(Y)] + R
    =R,
  \end{align*}
where
  \[
  |R| \le \frac{Md^{(k+1)/2}}{k!}|\eta|^{k+1},
  \]
and $M=\sup\{|\pa^\al g(u\eta)|:0\le u\le 1,|\al|=k+1\}$. Note that
  \[
  \pa^\al g(u\eta) = E[\pa^\al f(U + u\eta Y)Y^{|\al|}\check{\ph}_k(Y)]
    = E[\pa^\al f(\xi - \eta(1 - u)Y)Y^{|\al|}\check{\ph}_k(Y)].
  \]
Hence,
  \begin{align*}
  |\pa^\al g(u\eta)| &\le K\wt K
    E[(1 + |\xi - \eta(1 - u)Y|^r)|Y|^{|\al|}(1 + |Y|^r)]\\
  &\le K\wt K E[(1 + 2^r|\xi|^r + 2^r|\eta|^r|Y|^r)
    (|Y|^{|\al|}+|Y|^{|\al|+r}).
  \end{align*}
Since $|\eta|^2\le \nu d$, this completes the proof. \qed

\bigskip

The following special case will be used multiple times.

\begin{corollary}\label{C:Gauss_Taylor}
Let $X_1,\ldots,X_n$ be jointly normal, each with mean zero and variance bounded by $\nu>0$. Let $\eta_{ij}=E[X_iX_j]$. If $f\in C^1(\RR^{n-1})$ has polynomial growth of order 1 with constants $K$ and $r$, then
  \begin{equation}\label{Gauss_Taylor}
  |E[f(X_1,\ldots,X_{n-1})X_n^3]| \le C
  K
   \si^3\max_{j<n}|\eta_{jn}|,
  \end{equation}
where $\si=(EX_n^2)^{1/2}$ and $C$ depends only on
 $r$, $\nu$, and $n$.
\end{corollary}

\pf Apply Theorem \ref{T:Gauss_Taylor} with $k=0$. \qed

\bigskip

Finally, the following covariance estimates will be critical.

\begin{lemma}\label{L:covar}
Recall the notation $\be_j=(B(t_{j-1})+B(t_j))/2$ and $r_+=r\vee 1$. For any $i,j$,
  \begin{enumerate}[(i)]
  \item $|E[\De B_i\De B_j]| \le C\De t^{1/3}|j - i|_+^{-5/3}$,
  \item $|E[B(t_i)\De B_j]| \le C\De t^{1/3}(j^{-2/3} + |j - i|_+^{-2/3})$,
  \item $|E[\be_i\De B_j]|\le C\De t^{1/3}(j^{-2/3} + |j - i|_+^{-2/3})$,
  \item $|E[\be_j\De B_j]| \le C\De t^{1/3}j^{-2/3}$, and
  \item $C_1|t_j - t_i|^{1/3} \le E|\be_j - \be_i|^2
    \le C_2|t_j - t_i|^{1/3}$,
  \end{enumerate}
where $C_1,C_2$ are positive, finite constants that do not depend on $i$ or $j$.
\end{lemma}

\pf (i) By symmetry, we may assume $i\le j$. First, assume $j-i\ge 2$. Then
  \[
  E[\De B_i\De B_j] = \int_{t_{i-1}}^{t_i}\int_{t_{j-1}}^{t_j}
    \pa_{st}^2
    R(s,t)\,dt\,ds,
  \]
where $\pa^2_{st}=\pa_1\pa_2$. Note that for $s<t$, $\pa^2_{st}R(s,t) = -(1/9) (t - s)^{-5/3}$. Hence,
  \[
  |E[\De B_i\De B_j]| \le C\De t^2|t_{j-1} - t_i|^{-5/3}
    \le C\De t^{1/3}|j - i|^{-5/3}.
  \]
Now assume $j-i\le 1$. By H\"older's inequality, $|E[\De B_i\De B_j]| \le \De t^{1/3} = \De t^{1/3}|j - i|_+^{-5/3}$.

(ii) First note that by (i),
  \[
  |E[B(t_i)\De B_j]| \le \sum_{k=1}^i|E[\De B_k\De B_j]|
    \le C\De t^{1/3}\sum_{k=1}^j|k - j|_+^{-5/3}
    \le C\De t^{1/3}.
  \]
This
 proves the lemma when either $j=1$ or $|j-i|_+=1$. To complete the proof of (ii), suppose $j>1$ and $|j-i|>1$. Note that if $t>0$ and $s\ne t$, then
  \[
  \pa_2 R(s,t) = \frac16t^{-2/3} - \frac16|t - s|^{-2/3}\sgn(t - s).
  \]
We may therefore write $E[B(t_i)\De B_j] = \int_{t_{j-1}}^{t_j}\pa_2 R(t_i,u) \,du$, giving
  \[
  |E[B(t_i)\De B_j]|
    \le \De t\sup_{u\in[t_{j-1},t_j]}|\pa_2 R(t_i,u)|
    \le C\De t^{1/3}(j^{-2/3} + |j - i|_+^{-2/3}),
  \]
which is (ii).

(iii) This follows immediately from (ii).

(iv) Note that $2\be_j\De B_j=B(t_j)^2 -B(t_{j-1})^2$. Since $EB(t)^2=t^{1/3}$, the mean value theorem gives $|E [\be_j\De B_j]|\le C(\De t)t_j^{-2/3}=C\De t^{1/3}j^{-2/3}$.

(v) Without loss of generality, we may assume $i<j$. The upper bound follows from
  \[
  2(\be_j - \be_i) = (B(t_j) - B(t_i)) + (B(t_{j-1}) - B(t_{i-1})),
  \]
and the fact that $E|B(t)-B(s)|^2=|t-s|^{1/3}$. For the lower bound, we first assume $i<j-1$ and write
  \[
  2(\be_j - \be_i) = 2(B(t_{j-1}) - B(t_i)) + \De B_j + \De B_i.
  \]
Hence,
  \[
  (E|\be_j - \be_i|^2)^{1/2} \ge |t_{j-1} - t_i|^{1/6}
    - \frac12(E|\De B_j + \De B_i|^2)^{1/2}.
  \]
Since $\De B_i$ and $\De B_j$ are negatively correlated,
  \[
  E|\De B_j + \De B_i|^2 \le E|\De B_j|^2 + E|\De B_i|^2 = 2\De t^{1/3}.
  \]
Thus,
  \[
  (E|\be_j - \be_i|^2)^{1/2} \ge \De t^{1/6}|j - 1 - i|^{1/6}
    - 2^{-1/2}\De t^{1/6} \ge C\De t^{1/6}|j - i|^{1/6},
  \]
for some $C>0$. This completes the proof when $i<j-1$.

If $i=j-1$, the conclusion is immediate, since $2(\be_j - \be_{j-1}) = B(t_j) - B(t_{j-2})$. \qed

\subsection{Sextic and signed cubic variations}

\begin{thm}\label{T:sextic_var}
For each $T>0$, we have $E[\sup_{0\le t\le T}|V_n^6(B,t) - 15t|^2] \to 0$ as $n\to\infty$. In particular, $V_n^6(B,t)\to 15t$ ucp.
\end{thm}

\pf  Since $V_n^6(B)$ is monotone, it will suffice to show that $V_n^6(B,t) \to15t$ in $L^2$ for each fixed $t$. Indeed, the uniform convergence will then be a direct consequence of Dini's theorem. We write
  \[
  V_n^6(B,t) - 15t = \sum_{j=1}^{\flr{nt}}(\De B_j^6 - 15\De t)
    + 15(\flr{nt}/n - t).
  \]
Since $|\flr{nt}/n-t|\le\De t$, it will suffice to show that $E|\sum_{j=1} ^{\flr{nt}}(\De B_j^6 - 15\De t)|^2 \to 0$.
For this, we compute
  \begin{equation}\label{sextic_var}
  \begin{split}
  E\bigg|\sum_{j=1}^{\flr{nt}}(\De B_j^6 - 15\De t)\bigg|^2
    &= \sum_{i=1}^{\flr{nt}}\sum_{j=1}^{\flr{nt}}
    E[(\De B_i^6 - 15\De t)(\De B_j^6 - 15\De t)]\\
  &= \sum_{i=1}^{\flr{nt}}\sum_{j=1}^{\flr{nt}}
    (E[\De B_i^6\De B_j^6] - 225\De t^2).
  \end{split}
  \end{equation}
By Theorem \ref{T:Gauss_Taylor}, if $\xi,Y$ are jointly Gaussian, standard normals, then $E[\xi^6Y^6]=225+R$, where $|R|\le C|E[\xi Y]|^2$. Applying this with $\xi=\De t^{-1/6}\De B_i$ and $Y=\De t^{-1/6}\De B_j$, and using Lemma \ref{L:covar}(i), gives $|E[\De B_i^6\De B_j^6] - 225\De t^2]| \le C\De t^2|j - i|_+^{-10/3}$. Substituting this into \eqref{sextic_var}, we have
  \[
  E\bigg|\sum_{j=1}^{\flr{nt}}(\De B_j^6 - 15\De t)\bigg|^2
    \le C\flr{nt}\De t^2 \le Ct\De t \to 0,
  \]
which completes the proof. \qed

\begin{thm}\label{T:signed_cubic}
As $n\to\infty$, $(B,V_n(B))\to(B,\cub{B})$ in law in $D_{\RR^2}[0,\infty)$.
\end{thm}

\pf By Theorem 10 in \cite{NO}, $(B,V_n(B))\to(B,\kappa W)=(B,\cub{B})$
in law in $(D_\RR[0, \infty))^2$. By Lemma \ref{L:prod_space}, this implies $(B,V_n(B))\to(B,\cub{B})$ in $D_{\RR^2}[0,\infty)$. \qed

\subsection{Main result}\label{S:main}

Given $g\in C^\infty(\RR)$, choose $G$ such that $G'=g$. We then define
  \begin{equation}\label{main_int_def}
  \int_0^t g(B(s))\,dB(s) = G(B(t)) - G(B(0))
    + \frac1{12}\int_0^t G'''(B(s))\,d\cub{B}_s.
  \end{equation}
Note that, by definition, the change of variable formula \eqref{ito_form} holds for all $g\in C^\infty$. We shall use the shorthand notation $\int g(B)\,dB$ to refer to the process $t\mapsto\int_0^t g(B(s))\,dB(s)$. Similarly, $\int g(B)\,d\cub{B}$ and $\int g(B)\,ds$ shall refer to the processes $t\mapsto\int_0^t g(B(s))\,d\cub{B}_s$ and $t \mapsto \int_0^t g(B(s)) \,ds$, respectively.

Our main result is the following.

\begin{thm}\label{T:main}
If $g\in C^\infty(\RR)$, then $(B,V_n(B),I_n(g,B))\to(B,\cub{B}, \int g(B)\,dB)$ in law in $D_{\RR^3}[0,\infty)$.
\end{thm}

We also have the following generalization concerning the joint convergence of multiple sequences of Riemann sums.

\begin{thm}\label{T:joint_conv}
Fix $k\ge 1$. Let $g_j\in C^\infty(\RR)$ for $1\le j\le k$. Let $J_n$ be the $\RR^k$-valued process whose $j$-th component is $(J_n)_j =I_n(g_j,B)$. Similarly, define $J$ by $J_j=\int g_j(B)\,dB$. Then $(B,V_n(B),J_n)\to(B,\cub{B},J)$ in law in $D_{\RR^{k+2}}[0,\infty)$.
\end{thm}

\begin{remark}\label{R:joint_conv}
In less formal language, Theorem \ref{T:joint_conv} states that the Riemann sums $I_n(g_j,B)$ converge jointly, and the limiting stochastic integrals are all defined in terms of the same Brownian motion. In other words, the limiting Brownian motion remains unchanged under changes in the integrand. In this sense, the limiting Brownian motion depends only on $B$, despite being independent of $B$ in the probabilistic sense.
\end{remark}

The proofs of these two theorems are given in Section \ref{S:main_proof}.

\section{Finite-dimensional distributions}

\begin{thm}\label{T:main_fdd}
If $g\in C^\infty(\RR)$ is bounded with bounded derivatives, then
  \[
  \left(B,V_n(B),\frac1{\sqrt{n}}\sum_{j=1}^{\flr{n\cdot}}
    \frac{g(B(t_{j-1})) + g(B(t_j))}2
    h_3(n^{1/6}\De B_j)\right)
    \to \left(B,\cub{B},\int g(B)\,d\cub{B}\right),
  \]
in the sense of finite-dimensional distributions on $[0,\infty)$.
\end{thm}

The rest of this section is devoted to the proof of Theorem \ref{T:main_fdd}.

\subsection{Some technical lemmas}\label{sec:tec}

During the proof of Theorem \ref{T:main_fdd}, we will need technical results that are collected here. Moreover, for notational convenience, we will make use of the following shorthand notation:
  \[
  \de_j = {\bf 1}_{[t_{j-1},t_j]} \quad\text{ and }\quad
    \ep_j = {\bf 1}_{[0,t_{j}]}.
  \]
For future reference, let us note that by \eqref{symmetrization},
  \begin{equation}\label{sym-formula}
  \ep_t^{\ot a}\tot\ep_s^{\ot(q-a)} = \binom qa^{-1}
    \sum_{\substack{i_1,\ldots,i_q\in\{s,t\}\\
      |\{j:i_j = s\}| = q - a}}
    \ep_{i_1}\ot\ldots\ot\ep_{i_q}.
  \end{equation}

\begin{lemma}\label{lemma:Tec1}
We have
  \begin{enumerate}[(i)]
  \item $|E[B(r)(B(t) - B(s))]|
    = |\ang{{\bf 1}_{[0,r]},{\bf 1}_{[s,t]}}_\HH|
    \le |t - s|^{1/3}$ for any $r,s,t\ge 0$;
  \item $\ds{\sup_{0\le s\le T}\sum_{k=1}^{\flr{nt}}
    |E[B(s)\De B_k]| = \sup_{0\le s\le T}\sum_{k=1}^{\flr{nt}}
    |\ang{{\bf 1}_{[0,s]},\de_k}_\HH|
    \underset{n\to\infty}{=} O(1)}$ for any fixed $t,T>0$;
  \item $\ds{\sum_{k,j=1}^{\flr{nt}} |E(B(t_{j-1})\De B_k)|
    = \sum_{k,j=1}^{\flr{nt}} |\ang{\ep_{j-1},\de_k}_\HH|
    \underset{n\to\infty}{=} O(n)}$ for any fixed $t>0$;
  \item $\ds{\sum_{k=1}^{\flr{nt}}\left|{(E[B(t_{k-1})\De B_k])^3
    + \frac1{8n}}\right|
    = \sum_{k=1}^{\flr{nt}}\left|{\ang{\ep_{k-1},\de_k}_\HH^3
    + \frac1{8n}}\right|
    \underset{n\to\infty}{\longrightarrow} 0}$ for any fixed $t>0$;
  \item $\ds{\sum_{k=1}^{\flr{nt}}\left|{(E[B(t_k)\De B_k]\big)^3
    - \frac1{8n}}\right|
    = \sum_{k=1}^{\flr{nt}}\left|{\ang{\ep_{k},\de_k}_\HH^3
    - \frac1{8n}}\right|
    \underset{n\to\infty}{\longrightarrow} 0}$ for any fixed $t>0$.
  \end{enumerate}
\end{lemma}

\pf
\begin{enumerate}[$(i)$]
\item We have
  \[
  E\big(B(r)(B(t) - B(s))\big)
    = \frac12(t^{1/3} - s^{1/3})
    + \frac12\left(|s - r|^{1/3} - |t - r|^{1/3}\right).
  \]
Using the classical inequality $\big||b|^{1/3}-|a|^{1/3}\big|\le |b-a|^{1/3}$, the desired result follows.
\item Observe that
  \[
  E(B(s)\De B_k) = \frac1{2n^{1/3}}
    \left(k^{1/3} - (k - 1)^{1/3} - |k - ns|^{1/3}
    + |k - ns - 1|^{1/3} \right).
  \]
We deduce, for any fixed $s\le t$:
  \begin{multline*}
  \sum_{k=1}^{\flr{nt}} |E(B(s)\De B_k)| \le \frac12 t^{1/3}
    + \frac1{2n^{1/3}}\bigg(
    (\flr{ns} - ns + 1)^{1/3} - (ns - \flr{ns})^{1/3}\\
  + \sum_{k=1}^{\flr{ns}} ((ns + 1 - k)^{1/3} - (ns - k)^{1/3})
    + \sum_{k=\flr{ns}+2}^{\flr{nt}}((k - ns)^{1/3} - (k - ns - 1)^{1/3})
    \bigg)\\
  = \frac12(t^{1/3} + s^{1/3} + |t - s|^{1/3}) + R_n,
  \end{multline*}
where $|R_n|\le Cn^{-1/3}$, and $C$ does not depend on $s$ or $t$. The case where $s> t$ can be obtained similarly. Taking the supremum over $s\in[0,T]$ gives us $(ii)$.
\item is a direct consequence of $(ii)$.
\item We have
  \begin{multline*}
  \left|\big(E(B(t_{k-1})\De B_k)\big)^3 + \frac1{8n}\right|
    = \frac1{8n}\left(k^{1/3} - (k - 1)^{1/3}\right)\\
  \times\left|\left(k^{1/3} - (k - 1)^{1/3}\right)^2
    - 3(k^{1/3} - (k - 1)^{1/3})+3\right|.
  \end{multline*}
Thus, the desired convergence is immediately checked by combining the bound $0\le k^{1/3}-(k-1)^{1/3}\le 1$ with a telescoping sum argument.
\item The proof is very similar to the proof of $(iv)$. \qed
\end{enumerate}

\begin{lemma}\label{lemma:Tec2}
Let $s\ge 1$, and suppose that $\phi\in C^6(\RR^s)$ and $g_1,g_2\in C^6(\RR)$ have polynomial growth of order $6$, all with constants $K$ and $r$. Fix $a,b\in[0,T]$. Then
  \[
  \sup_{u_1,\ldots,u_s\in[0,T]}\sup_{n\ge 1}
    \sum_{i_1=1}^{\flr{na}}\sum_{i_2=1}^{\flr{nb}}
    \big|E\big(\phi(B(u_1),\ldots,B(u_s))
    g_1(B(t_{i_1-1}))g_2(B(t_{i_2-1}))
    I_3(\delta_{i_1}^{\ot 3})I_3(\delta_{i_2}^{\ot 3})\big)\big|
  \]
is finite.
\end{lemma}

\pf Let $C$ denote a constant depending only on $T$, $s$, $K$, and $r$, and whose value can change from one line to another. Define $f:\R^{s+3}\to\R$ by
  \[
  f(x) = \phi(x_1,\ldots,x_s)g_1(x_{s+1})g_2(x_{s+2})h_3(x_{s+3}).
  \]
Let $\xi_i=B(u_i)$, $i=1,\ldots,s$; $\xi_{s+1}=B(t_{i_1-1})$, $\xi_{s+2}=B(t_{i_2-1})$, $\xi_{s+3}=n^{1/6}\Delta B_{i_1}$, and $\eta_i=n^{1/6}E[\xi_i\Delta B_{i_2}]$. Applying Theorem \ref{T:Gauss_Taylor} with $k=5$, we obtain
  \begin{multline*}
  E\big(\phi(B(u_1),\ldots,B(u_s))
    g_1(B(t_{i_1-1}))g_2(B(t_{i_2-1}))
    I_3(\de_{i_1}^{\ot 3})I_3(\de_{i_2}^{\ot 3})\big)\\
  = \frac1n E\big(\phi(B(u_1),\ldots,B(u_s))
    g_1(B(t_{i_1-1}))g_2(B(t_{i_2-1}))
    h_3(n^{\frac16}\De B_{i_1})h_3(n^{\frac16}\De B_{i_2})\big)\\
  = \frac1n\sum_{|\al|= 3} \frac6{\al!}E[\pa^\al f(\xi)]\eta^\al
    + \frac Rn,
  \end{multline*}
where $|R|\leq C|\eta|^{6}$.

By Lemma \ref{lemma:Tec1} $(i)$, we have $|\eta_i|\le n^{-1/6}$ for any $i\le s+2$, and $|\eta_{s+3}|\le 1$. Moreover, we have
  \[
  \frac1n\sum_{i_1=1}^{\flr{na}}\sum_{i_2=2}^{\flr{nb}}|\eta_{s+3}|
    = \frac1{2n}\sum_{i_1=1}^{\flr{na}}\sum_{i_2=2}^{\flr{nb}}
    \big||i_1 - i_2+1|^{1/3} + |i_1 - i_2 - 1|^{1/3}
    - 2|i_1 - i_2|^{1/3}\big| \le C.
  \]
Therefore, by taking into account these two facts, we deduce $\frac1n \sum_{i_1=1}^{\flr{na}}\sum_{i_2=2}^{\flr{nb}}|R|\le C$.

On the other hand, if $\al\in\NN_0^{s+3}$ is such that $|\al|=3$ with $\al_{s+3}\ne 0$, we have
  \begin{multline*}
  \frac1n \sum_{i_1=1}^{\flr{na}}\sum_{i_2=1}^{\flr{nb}}
    \frac6{\al!}\big|E[\pa^\al f(\xi)]\big| |\eta^\al|\\
  \le \frac Cn
    \sum_{i_1=1}^{\flr{na}}\sum_{i_2=1}^{\flr{nb}}
    \big||i_1-i_2+1|^{1/3}+|i_1-i_2-1|^{1/3}-2|i_1-i_2|^{1/3}\big|
    \le C.
  \end{multline*}
Finally, if $\al\in\NN_0^{s+3}$ is such that $|\al|=3$ with $\al_{s+3}=0$ then $\pa^\al f = \pa^\al\wh f\ot h_3$ with $\wh f:\R^{s+2}\to\R$ defined by $\wh f(x)=\phi(x_1,\ldots,x_s)g_1(x_{s+1})g_2(x_{s+2})$. Hence, applying Theorem \ref{T:Gauss_Taylor} to $\wh f$ with $k=2$, we deduce, for $\wh\eta\in\NN_0^{s+2}$ defined by $\wh\eta_i=\eta_i$,
  \[
  \big|E[\pa^\al f(\xi)]\big| = \big|E[\pa^\al\wh f(\xi)
    h_3(n^{1/6}\De B_{i_1})]\big|
    \le C|\wh\eta|^3 \le Cn^{-1/2},
  \]
so that
  \[
  \frac1n\sum_{i_1=1}^{\flr{na}}\sum_{i_2=1}^{\flr{nb}}
    \frac6{\al!}\big|E[\pa^\al f(\xi)]\big||\eta^\al|
    = \frac1n\sum_{i_1=1}^{\flr{na}}\sum_{i_2=1}^{\flr{nb}}
    \frac6{\al!}\big|E[\pa^\al f(\xi)]\big||\wh\eta^\al|
    \le C.
  \]
The proof of Lemma \ref{lemma:Tec2} is done. \qed

\begin{lemma}\label{lm-lei-rule}
Let $g,h\in C^q(\RR)$, $q\geq 1$, having bounded derivatives, and fix $s,t\geq 0$. Set $\ep_t={\bf 1}_{[0,t]}$ and $\ep_s={\bf 1}_{[0,s]}$. Then $g(B(t))h(B(s))$ belongs in $\mathbb{D}^{q,2}$ and we have
  \begin{equation}\label{lei-rule}
  D^q\big(g(B(t))h(B(s))\big) = \sum_{a=0}^q\binom{q}{a}
    g^{(a)}(B(t))h^{(q-a)}(B(s))\,\ep_t^{\ot a}\tot\ep_s^{\ot (q-a)}.
  \end{equation}
\end{lemma}
\pf This follows immediately from \eqref{tensor-leibnitz}. \qed

\begin{lemma}\label{lemma:Tec3}
Fix an integer $r\geq 1$, and some real numbers $s_1,\ldots,s_r\geq 0$. Suppose $\ph\in C^\infty(\RR^r)$ and $g_j\in C^\infty(\RR)$, $j=1,2,3,4$, are bounded with bounded partial derivatives. For $i_1,i_2,i_3, i_4 \in \NN$, set $\Phi(i_1,i_2,i_3,i_4):=\ph(B_{s_1}, \ldots,B_{s_r}) \prod_{j=1}^4 g_j(B_{s_{i_j}})$. Then, for any fixed $a,b,c,d>0$, the following estimate is in order:
  \begin{equation}\label{chg}
  \sup_{n\ge 1}\sum_{i_1=1}^{\flr{na}}\sum_{i_2=1}^{\flr{nb}}
    \sum_{i_3=1}^{\flr{nc}}\sum_{i_4=1}^{\flr{nd}}
    \left|E\left(\Phi(i_1,i_2,i_3,i_4)
    I_3(\de_{i_1}^{\ot 3})I_3(\de_{i_2}^{\ot 3})
    I_3(\de_{i_3}^{\ot 3})I_3(\de_{i_4}^{\ot 3})
  \right)\right| < \infty.
  \end{equation}
\end{lemma}

\pf
Using the product formula \eqref{multiplication}, we have that $I_3(\de_{ i_3}^{\ot 3}) I_3(\de_{i_4}^{\ot 3})$ equals
  \[
  I_6(\de_{i_3}^{\ot 3} \ot \de_{i_4}^{\ot 3})
    + 9 I_4(\de_{i_3}^{\ot 2} \ot \de_{i_4}^{\ot 2})
      \ang{\de_{i_3},\de_{i_4}}_{\HH}
    + 18 I_2(\de_{i_3} \ot \de_{i_4})
      \ang{\de_{i_3},\de_{i_4}}_{\HH}^2
    + 6\ang{\de_{i_3},\de_{i_4}}_{\HH}^3.
  \]
As a consequence, we get
  \begin{align*}
  \sum_{i_1=1}^{\flr{na}}&\sum_{i_2=1}^{\flr{nb}}
    \sum_{i_3=1}^{\flr{nc}}\sum_{i_4=1}^{\flr{nd}}
    \left|E\left(\Phi(i_1,i_2,i_3,i_4)
    I_3(\de_{i_1}^{\ot 3})I_3(\de_{i_2}^{\ot 3})
    I_3(\de_{i_3}^{\ot 3})I_3(\de_{i_4}^{\ot 3})
    \right)\right|\\
  &\le \sum_{i_1=1}^{\flr{na}}\sum_{i_2=1}^{\flr{nb}}
    \sum_{i_3=1}^{\flr{nc}}\sum_{i_4=1}^{\flr{nd}}
    \left|E\left(\Phi(i_1,i_2,i_3,i_4)
    I_3(\de_{i_1}^{\ot 3})I_3(\de_{i_2}^{\ot 3})
    I_6(\de_{i_3}^{\ot 3}\ot\de_{i_4}^{\ot 3})
    \right)\right|\\
  &\quad + 9\sum_{i_1=1}^{\flr{na}}\sum_{i_2=1}^{\flr{nb}}
    \sum_{i_3=1}^{\flr{nc}}\sum_{i_4=1}^{\flr{nd}}
    \left|E\left(\Phi(i_1,i_2,i_3,i_4)
    I_3(\de_{i_1}^{\ot 3})I_3(\de_{i_2}^{\ot 3})
    I_4(\de_{i_3}^{\ot 2}\ot\de_{i_4}^{\ot 2})
    \right)\right|\left|\ang{\de_{i_3},\de_{i_4}}_{\HH}\right|\\
  &\quad + 18\sum_{i_1=1}^{\flr{na}}\sum_{i_2=1}^{\flr{nb}}
    \sum_{i_3=1}^{\flr{nc}}\sum_{i_4=1}^{\flr{nd}}
    \left|E\left(\Phi(i_1,i_2,i_3,i_4)
    I_3(\de_{i_1}^{\ot 3})I_3(\de_{i_2}^{\ot 3})
    I_2(\de_{i_3} \ot \de_{i_4})
    \right)\right|\ang{\de_{i_3},\de_{i_4}}_{\HH}^2\\
  &\quad + 6\sum_{i_1=1}^{\flr{na}}\sum_{i_2=1}^{\flr{nb}}
    \sum_{i_3=1}^{\flr{nc}}\sum_{i_4=1}^{\flr{nd}}
    \left|E\left(\Phi(i_1,i_2,i_3,i_4)
    I_3(\de_{i_1}^{\ot 3})I_3(\de_{i_2}^{\ot 3})
    \right)\right|\left|\ang{\de_{i_3},\de_{i_4}}_\HH\right|^3\\
  &=: A_1^{(n)} + 9A_2^{(n)} + 18A_3^{(n)} + 6A_4^{(n)}.
  \end{align*}

$(1)$ First, we deal with the term $A_1^{(n)}$.
  \begin{align*}
  A_1^{(n)} &= \sum_{i_1=1}^{\flr{na}}\sum_{i_2=1}^{\flr{nb}}
    \sum_{i_3=1}^{\flr{nc}}\sum_{i_4=1}^{\flr{nd}}
    \left\vert E\left( \Phi(i_1,i_2,i_3,i_4)
    I_3(\de_{i_1}^{\ot 3})I_3(\de_{i_2}^{\ot 3})
    I_6(\de_{i_3}^{\ot 3}\ot\de_{i_4}^{\ot 3})\right)\right\vert\\
  &= \sum_{i_1=1}^{\flr{na}}\sum_{i_2=1}^{\flr{nb}}
    \sum_{i_3=1}^{\flr{nc}}\sum_{i_4=1}^{\flr{nd}}
    \left\vert E\left(\ang{D^6\left(\Phi(i_1,i_2,i_3,i_4)
    I_3(\de_{i_1}^{\ot 3})I_3(\de_{i_2}^{\ot 3})\right),
    \de_{i_3}^{\ot 3}\ot\de_{i_4}^{\ot 3}}_{\HH^{\ot 6}}
    \right)\right\vert
  \end{align*}
When computing the sixth Malliavin derivative $D^6\big(\Phi(i_1,i_2,i_3,i_4) I_3(\de_{i_1}^{\ot 3}) I_3(\de_{i_2}^{\ot 3})\big)$, there are three types of terms:

$(1a)$ The first type consists in terms arising when one {\it only} differentiates $\Phi(i_1,i_2,i_3,i_4)$. By Lemma \ref{lemma:Tec1} $(i)$, these terms are all bounded by
  \[
  n^{-2}\sum_{i_1=1}^{\flr{na}}\sum_{i_2=1}^{\flr{nb}}
    \sum_{i_3=1}^{\flr{nc}}\sum_{i_4=1}^{\flr{nd}}
    \left\vert E\left(\wt\Phi(i_1,i_2,i_3,i_4)
    I_3(\de_{i_1}^{\ot 3})I_3(\de_{i_2}^{\ot 3})
    \right)\right\vert,
  \]
which is less than
  \[
  cd\sup_{i_3=1,\ldots,\flr{nc}}\sup_{i_4=1,\ldots,\flr{nd}}
    \sum_{i_1=1}^{\flr{na}}\sum_{i_2=1}^{\flr{nb}}
    \left\vert E\left(\wt\Phi(i_1,i_2,i_3,i_4)
    I_3(\de_{i_1}^{\ot 3}) I_3(\de_{i_2}^{\ot 3})
    \right)\right\vert.
  \]
(Here, $\wt\Phi(i_1,i_2,i_3,i_4)$ means a quantity having a similar form as $\Phi(i_1,i_2,i_3,i_4)$.) Therefore, Lemma \ref{lemma:Tec2} shows that the terms of the first type in $A^{(n)}_1$ well agree with the desired conclusion \eqref{chg}.

$(1b)$ The second type consists in terms arising when one differentiates $\Phi(i_1,i_2,i_3,i_4)$ and $I_3(\de_{i_1}^{\ot 3})$, but not $I_3( \de_{i_2}^{\ot 3})$ (the case where one differentiates $\Phi(i_1,i_2,i_3,i_4)$ and $I_3(\de_{i_2}^{\ot 3})$ but not $I_3(\de_{i_1}^{\ot 3})$ is, of course, completely similar). In this case, with $\rho$ defined by \eqref{rho}, the corresponding terms are bounded either by
  \[
  Cn^{-2}\sum_{i_1=1}^{\flr{na}}\sum_{i_2=1}^{\flr{nb}}
    \sum_{i_3=1}^{\flr{nc}}\sum_{i_4=1}^{\flr{nd}}
    \sum_{\al=0}^{2}\left\vert E\left(\wt\Phi(i_1,i_2,i_3,i_4)
    I_\al(\de_{i_1}^{\ot\al})I_3(\de_{i_2}^{\ot 3})
    \right)\right\vert|\rho(i_3 - i_1)|,
  \]
or by the same quantity with $\rho(i_4-i_1)$ instead of $\rho(i_3-i_1)$. In order to get the previous estimate, we have used Lemma \ref{lemma:Tec1} $(i)$ plus the fact that the sequence $\{\rho(r)\}_{r\in\ZZ}$, introduced in \eqref{rho}, is bounded. Moreover, by \eqref{dual} and Lemma \ref{lemma:Tec1}
$(i)$, observe that
  \[
  \left|E\left(\wt\Phi(i_1,i_2,i_3,i_4)
    I_\al(\de_{i_1}^{\ot\al})I_3(\de_{i_2}^{\ot 3})\right)\right|
    = \left|E\left(\ang{D^3\big(\wt\Phi(i_1,i_2,i_3,i_4)
    I_\al(\de_{i_1}^{\ot\al})\big),
    \de_{i_2}^{\ot 3}}_{\HH^{\otimes 3}}\right)\right|
    \le Cn^{-1},
  \]
for any $\al=0,1,2$. Finally, since
  \[
  \sup_{i_1=1,\ldots,\flr{na}}
    \sum_{i_3=1}^{\flr{nc}}\sum_{i_4=1}^{\flr{nd}}|\rho(i_3 - i_1)|
    \le nd\sup_{i_1=1,\ldots,\flr{na}}\sum_{r\in\ZZ}|\rho(r)|
    = Cn
  \]
(and similarly for $\rho(i_4-i_1)$ instead of $\rho(i_3-i_1)$), we deduce that the terms of the second type in $A^{(n)}_1$ also agree with the desired conclusion \eqref{chg}.

$(1c)$ The third and last type of terms consist of those that arise when one differentiates $\Phi(i_1,i_2,i_3,i_4)$, $I_3(\de_{i_1})$ and $I_3(\de_{i_2})$. In this case, the corresponding terms can be bounded by expressions of the type
  \[
  Cn^{-2}\sum_{i_1=1}^{\flr{na}}\sum_{i_2=1}^{\flr{nb}}
    \sum_{i_3=1}^{\flr{nc}}\sum_{i_4=1}^{\flr{nd}}
    \sum_{\al=0}^2\sum_{\be=0}^2
    \left\vert E\left(\wt\Phi(i_1,i_2,i_3,i_4)
    I_\al(\de_{i_1}^{\ot\al})I_\be(\de_{i_2}^{\ot\be})
    \right)\right\vert|\rho(i_3 - i_1)||\rho(i_2 - i_3)|.
  \]
Since $\left\vert E\left(\wt\Phi(i_1,i_2,i_3,i_4) I_\al(\de_{i_1}^{\ot \al})I_\be(\de_{i_2}^{\ot\be})\right)\right\vert$ is  uniformly bounded in $n$ on one hand, and
  \[
  \sum_{i_1=1}^{\flr{na}}\sum_{i_2=1}^{\flr{nb}}\sum_{i_3=1}^{\flr{nc}}
    |\rho(i_3 - i_1)||\rho(i_2-i_3)|
    \le nc\left(\sum_{r\in\ZZ}|\rho(r)|\right)^2 = Cn
  \]
on the other hand, we deduce that the terms of the third type in $A^{(n)}_1$ also agree with the desired conclusion \eqref{chg}.

$(2)$ Second, we focus on the term $A_2^{(n)}$. We have
  \[
  A_2^{(n)} = \sum_{i_1=1}^{\flr{na}}\sum_{i_2=1}^{\flr{nb}}
    \sum_{i_3=1}^{\flr{nc}}\sum_{i_4=1}^{\flr{nd}}
    \left\vert E\left(\ang{D^4\left(\Phi(i_1,i_2,i_3,i_4)
    I_3(\de_{i_1}^{\ot 3})I_3(\de_{i_2}^{\ot 3})\right),
    \de_{i_3}^{\ot 2}\ot\de_{i_4}^{\ot 2}}_{\HH^{\ot 4}}\right)\right|
    \left|\ang{\de_{i_3},\de_{i_4}}_\HH\right\vert.
  \]
When computing the fourth Malliavin derivative $D^4\big(\Phi(i_1,i_2,i_3,i_4) I_3(\de_{i_1}^{\ot 3}) I_3(\de_{i_2}^{\ot 3})\big)$, we have to deal with three types of terms:

$(2a)$ The first type consists in terms arising when one only differentiates $\Phi(i_1,i_2,i_3,i_4)$. By Lemma \ref{lemma:Tec1} $(i)$, these terms are all bounded by
  \[
  n^{-5/3}\sum_{i_1=1}^{\flr{na}}\sum_{i_2=1}^{\flr{nb}}
    \sum_{i_3=1}^{\flr{nc}}\sum_{i_4=1}^{\flr{nd}}
    \left\vert E\left(\wt\Phi(i_1,i_2,i_3,i_4)
    I_3(\de_{i_1}^{\ot 3})I_3(\de_{i_2}^{\ot 3})
    \right)\right\vert|\rho(i_3 - i_4)|,
  \]
which is less than
  \[
  Cn^{-2/3}\sum_{r\in\ZZ}|\rho(r)|
    \sup_{i_3=1,\ldots,\flr{nc}}\sup_{i_4=1,\ldots,\flr{nd}}
    \sum_{i_1=1}^{\flr{na}}\sum_{i_2=1}^{\flr{nb}}
    \left\vert E\left(\wt\Phi(i_1,i_2,i_3,i_4)
    I_3(\de_{i_1}^{\ot 3})I_3(\de_{i_2}^{\ot 3})
    \right)\right\vert.
  \]
Hence, by Lemma \ref{lemma:Tec2}, we see that the terms of the first type in $A^{(n)}_2$ well agree with the desired conclusion \eqref{chg}.

$(2b)$ The second type consists in terms arising when one differentiates $\Phi(i_1,i_2,i_3,i_4)$ and $I_3(\delta_{i_1}^{\otimes 3})$ but not $I_3(\delta_{i_2}^{\otimes 3})$ (the case where one differentiates $\Phi(i_1,i_2,i_3,i_4)$ and $I_3(\delta_{i_2}^{\otimes 3})$ but not $I_3(\delta_{i_1}^{\otimes 3})$ is completely similar). In this case, the corresponding terms can be bounded either by
  \[
  Cn^{-5/3}\sum_{i_1=1}^{\flr{na}}\sum_{i_2=1}^{\flr{nb}}
    \sum_{i_3=1}^{\flr{nc}}\sum_{i_4=1}^{\flr{nd}}\sum_{\al=0}^2
    \left\vert E\left(\wt\Phi(i_1,i_2,i_3,i_4)
    I_\al(\de_{i_1}^{\ot\al})I_3(\de_{i_2}^{\ot 3})
    \right)\right\vert|\rho(i_3 - i_1)||\rho(i_3 - i_4)|,
  \]
or by the same quantity with $\rho(i_4-i_1)$ instead of $\rho(i_3-i_1)$. By Cauchy-Schwarz inequality, we have
  \[
  \left\vert E\left(\wt\Phi(i_1,i_2,i_3,i_4)
    I_\al(\de_{i_1}^{\ot\al})I_3(\de_{i_2}^{\ot 3})\right)\right\vert
    \le Cn^{-\frac{3 + \al}6} \le Cn^{-1/2}.
  \]
Since moreover
  \[
  \sup_{i_1=1,\ldots,\flr{na}}
    \sum_{i_3=1}^{\flr{nc}}\sum_{i_4=1}^{\flr{nd}}
    |\rho(i_3 - i_1)| |\rho(i_3 - i_4)|
    \le \bigg(\sum_{r\in\ZZ} |\rho(r)|\bigg)^2 = C
  \]
(and similarly for $\rho(i_4-i_1)$ instead of $\rho(i_3-i_1)$), we deduce that the terms of the second type in $A^{(n)}_2$ also agree with the desired conclusion \eqref{chg}.

$(2c)$ The third and last type of terms consist of those that arise when one differentiates $\Phi(i_1,i_2,i_3,i_4)$, $I_3(\de_{i_1})$ and $I_3(\de_{i_2})$. In this case, the corresponding terms can be bounded by expressions of the type
  \begin{multline*}
  Cn^{-5/3}\sum_{i_1=1}^{\flr{na}}\sum_{i_2=1}^{\flr{nb}}
    \sum_{i_3=1}^{\flr{nc}}\sum_{i_4=1}^{\flr{nd}}
    \sum_{\al=0}^2\sum_{\be=0}^2
    \left\vert E\left(\wt\Phi(i_1,i_2,i_3,i_4)
    I_\al(\de_{i_1}^{\ot\al})I_\be(\de_{i_2}^{\ot\be})
    \right)\right\vert\\
  \times |\rho(i_3 - i_1)||\rho(i_2 - i_3)||\rho(i_3 - i_4)|.
  \end{multline*}
Since $\left\vert E\left(\wt\Phi(i_1,i_2,i_3,i_4) I_\al(\de_{i_1}^{\ot \al})I_\be(\de_{i_2}^{\ot\be})\right)\right\vert$ is uniformly bounded in $n$ on one hand, and
  \[
  \sum_{i_1=1}^{\flr{na}}\sum_{i_2=1}^{\flr{nb}}
    \sum_{i_3=1}^{\flr{nc}}\sum_{i_4=1}^{\flr{nd}}
    |\rho(i_3 - i_1)||\rho(i_2 - i_3)||\rho(i_3 - i_4)|
    \le nd\left(\sum_{r\in\ZZ}|\rho(r)|\right)^3 = Cn
  \]
on the other hand, we deduce that the terms of the third type in $A^{(n)}_2$ also agree with the desired conclusion \eqref{chg}.

$(3)$ Using exactly the same strategy than in point (2), we can show as well that the terms $A_3^{(n)}$ agree with the desired conclusion \eqref{chg}. Details are left to the reader.

$(4)$ Finally, let us focus on the last term, that is $A_4^{(n)}$. We have, using successively the fact that $\sum_{r\in\ZZ}|\rho(r)|^3 <\infty$ and Lemma \ref{lemma:Tec2},
  \begin{align*}
  A_4^{(n)} &= \sum_{i_1=1}^{\flr{na}}\sum_{i_2=1}^{\flr{nb}}
    \sum_{i_3=1}^{\flr{nc}}\sum_{i_4=1}^{\flr{nd}}
    \left\vert E\left(\Phi(i_1,i_2,i_3,i_4)
    I_3(\de_{i_1}^{\ot 3})I_3(\de_{i_2}^{\ot 3})\right)\right|
    \left|\ang{\de_{i_3},\de_{i_4}}_{\HH}\right\vert^3\\
  &= n^{-1}\sum_{i_1=1}^{\flr{na}}\sum_{i_2=1}^{\flr{nb}}
    \sum_{i_3=1}^{\flr{nc}}\sum_{i_4=1}^{\flr{nd}}
    \left\vert E\left(\Phi(i_1,i_2,i_3,i_4)
    I_3(\de_{i_1}^{\ot 3})I_3(\de_{i_2}^{\ot 3})\right)\right|
    |\rho(i_3-i_4)|^3\\
  &\le C \sup_{1\le i_3\le\flr{nc}}\sup_{1\le i_4\le\flr{nc}}
    \sum_{i_1=1}^{\flr{na}}\sum_{i_2=1}^{\flr{nb}}
    \left\vert E\left(\Phi(i_1,i_2,i_3,i_4)
    I_3(\de_{i_1}^{\ot 3})I_3(\de_{i_2}^{\ot 3})\right)\right\vert
    \le C.
  \end{align*}
Hence, the terms $A_4^{(n)}$ agree with the desired conclusion \eqref{chg} and the proof of Lemma \ref{lemma:Tec3} is now complete. \qed

\begin{lemma}\label{L:star}
Let $\lambda=(\lambda_1,\ldots,\lambda_m)\in\R^m$, $u_1,\ldots,u_m> 0$, $u_p> 0$ and suppose  $g_1,\ldots,g_m\in C^\infty(\RR)$ are bounded with bounded derivatives. Define $\VV_n\in\RR^m$ by
  \[
  \VV_n := \bigg(\sum_{i=1}^{\flr{nu_k}}
    g_k(B(t_{i-1}))I_3(\de_i^{\ot3})\bigg)_{k=1,\ldots,m},
  \]
so that
  \begin{equation}\label{star}
  \ang{\la,\VV_n} := \sum_{k=1}^m\la_k\sum_{i=1}^{\flr{nu_k}}
    g_k(B(t_{i-1}))I_3(\de_i^{\ot3})
    \quad\text{(see \eqref{vn} below)}.
  \end{equation}
Then there exists $C>0$, independent of $n$, such that
  \begin{align}
  \sup_{j=1,\ldots,\flr{nu_p}}
    E\big(\ang{D\ang{\la,\VV_n},\de_j}_\HH^2\big)
    &\le Cn^{-2/3}\label{est1}\\
  \sum_{j=1}^{\flr{nu_p}}
    E\big(\ang{D\ang{\la,\VV_n},\de_j}_\HH^2\big)
    &\le Cn^{-1/3}\label{est2}\\
  \sum_{j=1}^{\flr{nu_p}}
    E\big(\ang{D^2\ang{\la,\VV_n},\de_j^{\ot2}}_{\HH^{\ot2}}^2\big)
    &\le Cn^{-2/3}.\label{est3}
  \end{align}
\end{lemma}

\pf We have
  \begin{multline}\label{eq:DLambdaVn}
  \ang{D\ang{\la,\VV_n},\de_j}_\HH
    = \sum_{k=1}^m \la_k \sum_{i=1}^{\flr{nu_k}}
    g'_k(B(t_{i-1}))I_3(\de_i^{\ot3})\ang{\ep_{i-1},\de_j}_\HH\\
  + 3\sum_{k=1}^m \la_k \sum_{i=1}^{\flr{nu_k}}
    g_k(B(t_{i-1}))I_2(\de_i^{\ot2})\ang{\de_i,\de_j}_\HH.
  \end{multline}
Hence, with $\rho$ defined by \eqref{rho},
  \begin{align*}
  E&\big(\ang{D\ang{\la,\VV_n},\de_j}_\HH^2\big)\\
  &\le 2m\sum_{k=1}^m\la_k^2\sum_{i,\ell=1}^{\flr{nu_k}}
    \big|E(g'_k(B(t_{i-1}))g'_\ell(B(t_{\ell-1}))
    I_3(\de_i^{\ot3})I_3(\de_l^{\ot3}))\big|
    \big|\ang{\ep_{i-1},\de_j}_\HH\big|
    \big|\ang{\ep_{\ell-1},\de_j}_\HH\big|\\
  &\quad + 18m\sum_{k=1}^m\la_k^2\sum_{i,\ell=1}^{\flr{nu_k}}
    \big|E(g_k(B(t_{i-1}))g_k(B(t_{\ell-1}))
    I_2(\de_i^{\ot2})I_2(\de_l^{\ot2}))\big|
    \big|\ang{\de_i,\de_j}_\HH\big|
    \big|\ang{\de_\ell,\de_j}_\HH\big|\\
  &\le Cn^{-2/3}\sup_{k=1,\ldots,m}\sum_{i,\ell=1}^{\flr{nu_k}}
    \big|E(g'_k(B(t_{i-1}))g'_\ell(B(t_{\ell-1}))
    I_3(\de_i^{\ot3})I_3(\de_\ell^{\ot3}))\big|\\
  &\quad + Cn^{-4/3}\sum_{i,\ell=1}^{\flr{nu_k}}
    |\rho(i - j)||\rho(\ell - j)|
    \quad\text{by Lemma \ref{lemma:Tec1} $(i)$ and Cauchy-Schwarz}\\
  &\le Cn^{-2/3} + Cn^{-4/3}\bigg(\sum_{r\in\ZZ}|\rho(r)|\bigg)^2
    \quad\text{by Lemma \ref{lemma:Tec2}}\\
  &\le Cn^{-2/3},
  \end{align*}
which is \eqref{est1}. Moreover, combining the first inequality of the previous estimate with Lemma \ref{lemma:Tec1} $(ii)$ and Lemma \ref{lemma:Tec2}, we also have
  \begin{multline*}
  \sum_{j=1}^{\flr{nu_p}}E\big(\ang{D\ang{\la,\VV_n},\de_j}_\HH^2\big)\\
  \le Cn^{-1/3}\sup_{k=1,\ldots,m}\sum_{i,\ell=1}^{\flr{nu_k}}
    \big|E(g'_k(B(t_{i-1}))g'_\ell(B(t_{\ell-1}))
    I_3(\de_i^{\ot3})I_3(\de_l^{\ot3}))\big|\\
  \times \sup_{i=1,\ldots,\flr{nu_k}}\sum_{j=1}^{\flr{nu_p}}
    \big|\ang{\ep_{i-1},\de_j}_\HH\big|
    + Cn^{-1/3}\bigg(\sum_{r\in\ZZ}|\rho(r)|\bigg)^2
    \le Cn^{-1/3},
  \end{multline*}
which is \eqref{est2}. The proof of \eqref{est3} follows the same lines, and is left to the reader.\qed

\subsection{Proof of Theorem \ref{T:main_fdd}}
\label{section:mainfdd}

We are now in position to prove Theorem \ref{T:main_fdd}. For $g:\R \to\R$, let
  \[
  G_n^-(g,B,t) := \frac1{\sqrt{n}}\sum_{j=1}^{\flr{nt}}
    g(B(t_{j-1}))h_3(n^{1/6}\De B_j),
    \quad t\ge 0, \quad n\ge 1.
  \]
We recall that $h_3(x)=x^3-3x$, see \eqref{hermite}, and the definition (\ref{vn3}) of $V_n(B,t)$. In particular, observe that
  \begin{equation}\label{yahoo}
  V_n(B,t) = G_n^-(1,B,t) + 3 n^{-1/3} B(\flr{nt}/n).
  \end{equation}
Our main theorem which will lead us toward the proof of Theorem \ref{T:main_fdd} is the following.

\begin{thm}\label{fdd2}
If $g\in C^\infty(\RR)$ is bounded with bounded derivatives, then the sequence $(B,G_n^-(1,B),G_n^-(g,B))$ converges to $(B,\cub{B}, -(1/8)\int g'''(B)\,ds+\int g(B)\,d\cub{B})$ in the sense of finite-dimensional distributions on $[0,\infty)$.
\end{thm}

\pf We have to prove that, for any $\ell+m\ge 1$ and any $u_1,\ldots, u_{\ell+m}\ge 0$:
  \begin{multline*}
  \big(B,G_n^-(1,B,u_1),\ldots,G_n^-(1,B,u_\ell),
    G_n^-(g,B,u_{\ell+1}),\ldots,G_n^-(g,B,u_{\ell+m})\big)\\
  \xrightarrow[n\to\infty]{\text{Law}}
    \bigg(B,\cub{B}_{u_1},\ldots,\cub{B}_{u_\ell},
    - \frac18\int_0^{u_{\ell+1}}g'''(B(s))\,ds
    + \int_0^{u_{\ell+1}}g(B(s))\,d\cub{B}_s,\ldots,\\
  - \frac18\int_0^{u_{\ell+m}}g'''(B(s))\,ds
    + \int_0^{u_{\ell+m}}g(B(s))\,d\cub{B}_s\bigg).
  \end{multline*}
Actually, we will prove the following slightly stronger convergence. For any $m\ge 1$, any $u_1,\ldots,u_m\ge 0$ and all bounded functions $g_1, \ldots,g_m\in C^\infty(\RR)$ with bounded derivatives, we have
  \begin{multline}\label{afaire}
  \big(B,G_n^-(g_1,B,u_1),\ldots,G_n^-(g_m,B,u_m)\big)\\
  \xrightarrow[n\to\infty]{\text{Law}}
    \bigg(B,-\frac18\int_0^{u_1}g_1'''(B(s))\,ds
    + \int_0^{u_1}g_1(B(s))\,d\cub{B}_s,\ldots,\\
  - \frac18\int_0^{u_m}g_m'''(B(s))\,ds
    + \int_0^{u_m}g_m(B(s))\,d\cub{B}_s\bigg).
  \end{multline}
Using \eqref{eq:defIn}, observe that
  \begin{equation}\label{newform}
  G_n^-(g,B,t) = \sum_{j=1}^{\flr{nt}}
    g(B(t_{j-1}))I_3(\de_j^{\ot3}).
  \end{equation}
The proof of \eqref{afaire} is divided into several steps, and follows the methodology introduced in \cite{NN}.

\bigskip

\textit{Step 1.- } We first prove that:
  \begin{equation}\label{expectation}
  \begin{split}
  \lim_{n\to\infty} &E\left(
    G_n^-(g_1,B,u_1),\ldots,G_n^-(g_m,B,u_m)\right)\\
  &= \left(-\frac18 \int_0^{u_1} E(g_1'''(B(s)))\,ds,\ldots,
    -\frac18 \int_0^{u_m} E(g_m'''(B(s)))\,ds\right),\\
  \\
  \lim_{n\to\infty} &E\left(\left\|\big(
    G_n^-(g_1,B,u_1),\ldots,G_n^-(g_m,B,u_m)
    \big)\right\|_{\RR^m}^2\right)\\
  &= \sum_{i=1}^m \left(
    \ka^2 \int_0^{u_i} E(g_i^2(B(s)))\,ds
    + \frac1{64} E\left(\int_0^{u_i} g_i'''(B(s))\,ds\right)^2
    \right).
  \end{split}
  \end{equation}
For $g$ as in the statement of the theorem, we can write, for any fixed $t\ge 0$:
  \begin{align*}
  E\left(G_n^-(g,B,t)\right) &= \sum_{j=1}^{\flr{nt}}
    E\left(g(B(t_{j-1}))I_3(\de_j^{\ot3})\right)
    \quad\text{by \eqref{newform}}\\
  &= \sum_{j=1}^{\flr{nt}} E\left(
    \ang{D^3 g(B(t_{j-1})),\de_j^{\ot3}}_{\HH^{\ot3}}
    \right)\quad\text{by \eqref{dual}}\\
  &= \sum_{j=1}^{\flr{nt}}
    E(g'''(B(t_{j-1})))\ang{\ep_{j-1},\de_j}_\HH^3
    \quad\text{by \eqref{chainrule}}\\
  &= -\frac1{8n}\sum_{j=1}^{\flr{nt}} E(g'''(B(t_{j-1})))
    + \sum_{j=1}^{\flr{nt}} E(g'''(B(t_{j-1})))
    \left(\ang{\ep_{j-1},\de_j}_\HH^3 + \frac1{8n}\right)\\
  &\xrightarrow[n\to\infty]{} -\frac18\int_0^t E(g'''(B(s)))\,ds
    \quad\text{by Lemma \ref{lemma:Tec1} $(iv)$}.
  \end{align*}
Now, let us turn to the second part of \eqref{expectation}. We have
  \[
  E\|(G_n^-(g_1,B,u_1),\ldots,G_n^-(g_m,B,u_m))\|_{\RR^m}^2
    = \sum_{i=1}^m E(G_n^-(g_i,B,u_i)^2).
  \]
By the product formula \eqref{multiplication}, we have
  \begin{multline*}
  I_3(\de_j^{\ot3})I_3(\de_k^{\ot3})
    = I_6(\de_j^{\ot3}\ot\de_k^{\ot3})
    + 9I_4(\de_j^{\ot2}\ot\de_k^{\ot2})\ang{\de_j,\de_k}_\HH\\
  + 18I_2(\de_j\ot\de_k)\ang{\de_j,\de_k}_\HH^2
  + 6\ang{\de_j,\de_k}_\HH^3.
  \end{multline*}
Thus, for any fixed $t\ge 0$,
  \begin{align*}
  E(G_n^-(g,B,t)^2) &= \sum_{j,k=1}^{\flr{nt}} E\left(
    g(B(t_{j-1}))g(B(t_{k-1}))
    I_3(\de_j^{\ot3})I_3(\de_k^{\ot3})\right)\\
  &= \sum_{j,k=1}^{\flr{nt}} E\left(
    g(B(t_{j-1}))g(B(t_{k-1}))
    I_6(\de_j^{\ot3}\ot\de_k^{\ot3})\right)\\
  &\quad + 9\sum_{j,k=1}^{\flr{nt}} E\left(
    g(B(t_{j-1}))g(B(t_{k-1}))
    I_4(\de_j^{\ot2}\ot\de_k^{\ot2})\right)\ang{\de_j,\de_k}_\HH\\
  &\quad + 18\sum_{j,k=1}^{\flr{nt}} E\left(
    g(B(t_{j-1}))g(B(t_{k-1}))
    I_2(\de_j\ot\de_k)\right)\ang{\de_j,\de_k}_\HH^2\\
  &\quad + 6\sum_{j,k=1}^{\flr{nt}} E\left(
    g(B(t_{j-1}))g(B(t_{k-1}))\right)\ang{\de_j,\de_k}_\HH^3\\
  &=: A_n + B_n + C_n + D_n.
  \end{align*}
We will estimate each of these four terms using the Malliavin integration by parts formula \eqref{dual}. For that purpose, we use Lemma \ref{lm-lei-rule} and the notation of Remark \ref{rem-lei-rule}.

First, we have
  \begin{align*}
  A_n &= \sum_{j,k=1}^{\flr{nt}} E\left(\ang{
    D^6[g(B(t_{j-1}))g(B(t_{k-1}))],\de_j^{\ot3}\ot\de_k^{\ot3}
    }_{\HH^{\ot6}}\right)\\
  &\overset{\text{\eqref{lei-rule}}}{=}
    \sum_{j,k=1}^{\flr{nt}}\sum_{a=0}^6\binom6aE\left(
    g^{(a)}(B(t_{j-1}))g^{(6-a)}(B(t_{k-1}))\right)
    \ang{\ep_{j-1}^{\ot a}\tot\ep_{k-1}^{\ot(6-a)},
    \de_j^{\ot3}\ot\de_k^{\ot3}}_{\HH^{\ot6}}\\
  &\overset{\text{\eqref{sym-formula}}}{=}
    \sum_{j,k=1}^{\flr{nt}}\sum_{a=0}^6 E\left(
    g^{(a)}(B(t_{j-1}))g^{(6-a)}(B(t_{k-1}))\right)\\
  &\qquad\qquad\qquad\times
    \sum_{\substack{i_1,\ldots,i_6\in\{j-1,k-1\}\\
      |\{\ell: i_\ell = j - 1\}| = a}}
    \ang{\ep_{i_1}\ot\ldots\ot\ep_{i_6},
    \de_j^{\ot3}\ot\de_k^{\ot3}}_{\HH^{\ot6}}.
  \end{align*}
Actually, in the previous double sum with respect to $a$ and $i_1, \ldots,i_6$, only  the following term is non-negligible:
  \begin{align*}
  \sum_{j,k=1}^{\flr{nt}}E(g'''(B(&t_{j-1}))g'''(B(t_{k-1})))
    \ang{\ep_{j-1},\de_j}_\HH^3\ang{\ep_{k-1},\de_k}_\HH^3\\
  &= E\bigg(\sum_{j=1}^{\flr{nt}}
    g'''(B(t_{j-1}))\ang{\ep_{j-1},\de_j}_\HH^3\bigg)^2\\
  &= E\bigg(-\frac1{8n}\sum_{j=1}^{\flr{nt}}g'''(B(t_{j-1}))
    + \sum_{j=1}^{\flr{nt}}g'''(B(t_{j-1}))
    \bigg(\ang{\ep_{j-1},\de_j}_\HH^3 + \frac1{8 n}\bigg)\bigg)^2\\
  &\xrightarrow[n\to\infty]{} \frac1{64} E\bigg(
    \int_0^t g'''(B(s))\,ds\bigg)^2
    \quad\text{by Lemma \ref{lemma:Tec1} $(iv)$}.
  \end{align*}
Indeed, the other terms in $A_n$ are all of the form
  \begin{equation}\label{qu}
  \sum_{j,k=1}^{\flr{nt}}
    E(g^{(a)}(B(t_{j-1}))g^{(6-a)}(B(t_{k-1})))
    \ang{\ep_{j-1},\de_k}_\HH
    \prod_{i=1}^5\ang{\ep_{x_i-1},\de_{y_i}}_\HH,
  \end{equation}
where $x_i$ and $y_i$ are for $j$ or $k$. By Lemma \ref{lemma:Tec1} \textit{(iii)}, we have $\sum_{j,k=1}^{\flr{nt}}|\ang{\ep_{j-1}, \de_k}_\HH|=O(n)$ as $n\to\infty$. By Lemma \ref{lemma:Tec1} \textit{(i)},
$\sup_{j,k=1,\ldots,[nt]}\prod_{i=1}^5|\ang{\ep_{x_i-1},\de_{y_i} }_\HH|=O(n^{-5/3})$ as $n\to\infty$. Hence, the quantity in \eqref{qu} tends to zero as $n\to\infty$. We have proved
  \[
  A_n \xrightarrow[n\to\infty]{}
    \frac1{64} E\bigg(\int_0^t g'''(B(s))\,ds\bigg)^2.
  \]

Using the integration by parts formula \eqref{dual} as well as Lemma \ref{lm-lei-rule}, we have similarly that
  \begin{align*}
  |B_n| &\le \sum_{j,k=1}^{\flr{nt}}\sum_{a=0}^4
    \binom4a|E(g^{(a)}(B(t_{j-1}))g^{(4-a)}(B(t_{k-1})))
    \ang{\ep_{j-1}^{\ot a}\tot\ep_{k-1}^{\ot(4-a)},
    \de_j^{\ot2}\ot\de_k^{\ot2}}_{\HH^{\ot4}}
    \ang{\de_j,\de_k}_\HH|\\
  &\le C n^{-4/3}\sum_{j,k=1}^{\flr{nt}}
    |\ang{\de_j,\de_k}_\HH|
    \quad \text{by Lemma \ref{lemma:Tec1} $(i)$}\\
  &= C n^{-5/3}\sum_{j,k=1}^{\flr{nt}} |\rho(j-k)|
    \le C n^{-2/3} \sum_{r\in\ZZ} |\rho(r)|
    = Cn^{-2/3}
    \xrightarrow[n\to\infty]{} 0,
  \end{align*}
with $\rho$ defined by \eqref{rho}.

Using similar computations, we also have
  \[
  |C_n| \le Cn^{-1/3}\sum_{r=-\infty}^\infty\rho^2(r)
    = Cn^{-1/3} \xrightarrow[n\to\infty]{} 0,
  \]
while
  \begin{align*}
  D_n &= \frac6n\sum_{j,k=1}^{\flr{nt}}
    E(g(B(t_{j-1}))g(B(t_{k-1})))\rho^3(j-k)\\
  &= \frac6n\sum_{r\in\ZZ}
    \sum_{j=1\vee(1-r)}^{\flr{nt}\wedge(\flr{nt}-r)}
    E(g(B(t_{j-1}),)g(B(t_{j+r-1})))\rho^3(r)\\
  &\xrightarrow[n\to\infty]{} 6\sum_{r\in\ZZ}\rho^3(r)
    \int_0^t E(g^2(B(s)))\,ds
    = \kappa^2\int_0^t E(g^2(B(s)))\,ds,
  \end{align*}
the previous convergence being obtained as in the proof of \eqref{riemann2} below. Finally, we have obtained
  \begin{equation}\label{EGn2}
  E\big(G_n^-(g,B,t)^2\big)
    \xrightarrow[n\to\infty]{} \ka^2\int_0^t E\big(g^2(B(s))\big)\,ds     + \frac1{64} E\bigg(\int_0^t g'''(B(s))\,ds\bigg)^2,
  \end{equation}
and the proof of \eqref{expectation} is done.

\bigskip

{\it Step 2.-}
By Step 1, the sequence $\big(B,G_n^-(g_1,B,u_1),\ldots,G_n^-(g_m,B,u_m) \big)$ is tight in $D_\RR[0,\infty)\times\RR^m$. Consider a subsequence converging in law to some limit denoted by
  \[
  \big(B,G_\infty^-(g_1,B,u_1),\ldots,G_\infty^-(g_m,B,u_m)\big)
  \]
(for convenience, we keep the same notation for this subsequence and for the sequence itself). Recall $\VV_n$, defined in Lemma \ref{L:star}, and note that by \eqref{newform}, we have
  \begin{equation}\label{vn}
  \VV_n := \big(G_n^-(g_1,B,u_1),\ldots,G_n^-(g_m,B,u_m)\big),
    \quad n\in\NN\cup\{\infty\}.
  \end{equation}
Let us also define
  \begin{multline*}
  \WW := \bigg(-\frac18\int_0^{u_1} g_1'''(B(s))\,ds
    + \int_0^{u_1} g_1(B(s))\,d\cub{B}_s,\ldots,\\
  -\frac18\int_0^{u_m} g'''_m(B(s))\,ds
    + \int_0^{u_m} g_m(B(s))\,d\cub{B}_s\bigg).
  \end{multline*}
We have to show that, conditioned on $B$, the laws of $\VV_\infty$ and $\WW$ are the same.

Let $\la=(\la_1,\ldots,\la_m)$ denote a generic element of $\RR^m$ and, for $\la,\mu\in\RR^m$, write $\ang{\la,\mu}$ for $\sum_{i=1}^m \la_i\mu_i$. We consider the conditional characteristic function of $\WW$ given $B$:
  \begin{equation}\label{eq:definitioncaracteristicsfunctions}
  \Phi(\la) := E\left(e^{i\ang{\la,\WW}}\big\vert B\right).
  \end{equation}
Observe that $\Phi(\la)=e^{i\ang{\la,\mu}-\frac12\ang{\la,Q\la}}$, where $\mu_k:=-(1/8)\int_0^{u_k} g_k'''(B(s))\,ds$ for $k=1,\ldots,m$, and $Q=(q_{ij})_{1\le i,j\le m}$ is the symmetric matrix given by
  \[
  q_{ij} := \ka^2\int_0^{u_i\wedge u_j} g_i(B(s))g_j(B(s))\,ds.
  \]
The point is that $\Phi$ is the unique solution of the following system of PDEs (see \cite{R2}):
  \begin{equation}\label{eq:s^m}
  \frac{\pa\ph}{\pa\la_p}(\la)
    = \ph(\la)\bigg(i\mu_p - \sum_{k=1}^m\la_k q_{pk}\bigg),
    \quad p=1,\ldots,m,
  \end{equation}
where the unknown function $\ph:\RR^m\to\CC$ satisfies the initial condition $\ph(0)=1$. Hence, we have to show that, for every random variable $\xi$ of the form $\psi(B(s_1),\ldots,B(s_r))$, with $\psi: \RR^r \to \RR$ belonging to $C_b^\infty(\RR^r)$ and $s_1,\ldots,s_r \ge 0$, we have
  \begin{multline}\label{todo}
  \frac{\pa}{\pa\la_p} E\left(e^{i\ang{\la,\VV_\infty}}\xi\right)
    = -\frac i8\int_0^{u_p} E\left(g_p'''(B(s))
    \xi e^{i\ang{\la,\VV_\infty}}\right)\,ds\\
  - \ka^2\sum_{k=1}^m \la_k \int_0^{u_p\wedge u_k}
    E\left(g_p(B(s))g_k(B(s))
   \xi e^{i\ang{\la,\VV_\infty}}\right)\,ds
  \end{multline}
for all $p\in\{1,\ldots,m\}$.

\bigskip

{\it Step 3.-} Since $(\VV_\infty,B)$ is defined as the limit in law of $(\VV_n,B)$ on one hand, and $\VV_n$ is bounded in $L^2$ on the other hand, note that
  \[
  \frac{\pa}{\pa\la_p} E\left(e^{i\ang{\la,\VV_\infty}}\xi\right)
    = \lim_{n\to\infty} \frac{\pa}{\pa\la_p}
    E\left(e^{i\ang{\la,\VV_n}}\xi\right).
  \]
Let us compute $\frac{\pa}{\pa\la_p}E\left(e^{i\ang{\la,\VV_n}}\xi \right)$. We have
  \begin{equation}\label{eq:computationPartialPhi1}
  \frac{\pa}{\pa\la_p}E\left(e^{i\ang{\la,\VV_n}}\xi\right)
    = iE\big(G_n^-(g_p,B,u_p)e^{i\ang{\la,\VV_n}}\xi\big).
  \end{equation}
Moreover, see \eqref{newform} and use \eqref{dual}, for any $t\ge 0$:
  \begin{equation}\label{eq:computationPartialPhi}
  \begin{split}
  E\big(G_n^-(g,B,t)e^{i\ang{\la,\VV_n}}\xi\big)
    &= \sum_{j=1}^{\flr{nt}} E\left(g(B(t_{j-1}))
    I_3(\de_j^{\ot3})e^{i\ang{\la,\VV_n}}\xi\right)\\
  &= \sum_{j=1}^{\flr{nt}} E\left(\ang{D^3\left(
    g(B(t_{j-1}))e^{i\ang{\la,\VV_n}}\xi
    \right),\de_j^{\ot3}}_{\HH^{\ot3}}\right).
  \end{split}
  \end{equation}
The first three Malliavin derivatives of $g(B(t_{j-1}))e^{i\ang{\la, \VV_n}}\xi$ are respectively given by
  \begin{align*}
  D(g(&B(t_{j-1}))e^{i\ang{\la,\VV_n}}\xi)\\
  &= g'(B(t_{j-1}))e^{i\ang{\la,\VV_n}}\xi\,\,\ep_{j-1}
    + i g(B(t_{j-1}))e^{i\ang{\la,\VV_n}}\xi D\ang{\la,\VV_n}\\
  &\quad + g(B(t_{j-1}))e^{i\ang{\la,\VV_n}}D\xi,\\
  \\
  D^2(g(&B(t_{j-1}))e^{i\ang{\la,\VV_n}}\xi)\\
  &= g''(B(t_{j-1}))e^{i\ang{\la,\VV_n}}\xi\;\ep_{j-1}^{\ot2}
    + 2i g'(B(t_{j-1}))e^{i\ang{\la,\VV_n}}\xi
    \;D\ang{\la,\VV_n}\tot\ep_{j-1}\\
  &\quad + 2 g'(B(t_{j-1}))e^{i\ang{\la,\VV_n}}\;D\xi\tot\ep_{j-1}
    - g(B(t_{j-1}))e^{i\ang{\la,\VV_n}}\xi\;D\ang{\la,\VV_n}^{\ot2}\\
  &\quad + 2i g(B(t_{j-1}))e^{i\ang{\la,\VV_n}}
    \;D\xi\tot D\ang{\la,\VV_n}
    + i g(B(t_{j-1}))e^{i\ang{\la,\VV_n}}\xi\;D^2\ang{\la,\VV_n}\\
  &\quad + g(B(t_{j-1}))e^{i\ang{\la,\VV_n}}\;D^2\xi,
  \end{align*}
and
  \begin{equation}\label{eq:computationD3}
  \begin{split}
  D^3(g(&B(t_{j-1}))e^{i\ang{\la,\VV_n}}\xi)\\
  &= g'''(B(t_{j-1}))e^{i\ang{\la,\VV_n}}\xi\;\ep_{j-1}^{\ot3}
    + 3i g''(B(t_{j-1}))e^{i\ang{\la,\VV_n}}\xi
    \;\ep_{j-1}^{\ot2}\tot D\ang{\la,\VV_n}\\
  &\quad + 3 g''(B(t_{j-1}))e^{i\ang{\la,\VV_n}}
    \;\ep_{j-1}^{\ot2}\tot D\xi
    - 3 g'(B(t_{j-1}))e^{i\ang{\la,\VV_n}}\xi
    \;D\ang{\la,\VV_n}^{\ot2}\tot\ep_{j-1}\\
  &\quad + 6i g'(B(t_{j-1}))e^{i\ang{\la,\VV_n}}
    \;D\xi\tot D\ang{\la,\VV_n}\tot\ep_{j-1}\\
  &\quad - i g(B(t_{j-1}))e^{i\ang{\la,\VV_n}}\xi
    \;D\ang{\la,\VV_n}^{\ot3}
    - 3 g(B(t_{j-1}))e^{i\ang{\la,\VV_n}}
    \;D\ang{\la,\VV_n}^{\ot2}\tot D\xi\\
  &\quad + i g(B(t_{j-1}))e^{i\ang{\la,\VV_n}}\xi\;D^3\ang{\la,\VV_n}
    + g(B(t_{j-1}))e^{i\ang{\la,\VV_n}}\;D^3\xi\\
  &\quad + 3i g(B(t_{j-1}))e^{i\ang{\la,\VV_n}}
    \;D^2\xi\tot D\ang{\la,\VV_n}
    + 3i g(B(t_{j-1}))e^{i\ang{\la,\VV_n}}
    \;D\xi\tot D^2\ang{\la,\VV_n}\\
  &\quad + 3 g'(B(t_{j-1}))e^{i\ang{\la,\VV_n}}\;D^2\xi\tot\ep_{j-1}
    + 3i g'(B(t_{j-1}))e^{i\ang{\la,\VV_n}}\xi
    \;\ep_{j-1}\tot D^2\ang{\la,\VV_n}\\
  &\quad - 3 g(B(t_{j-1}))e^{i\ang{\la,\VV_n}}\xi\;
    D\ang{\la,\VV_n}\tot D^2\ang{\la,\VV_n}.
  \end{split}
  \end{equation}
Let us compute the term $D^3\ang{\la,\VV_n}$. Recall that
  \[
  \ang{\la,\VV_n} = \sum_{k=1}^m\la_k\,G_n^-(g_k,B,u_k)
    = \sum_{k=1}^m\la_k\sum_{\ell=1}^{\flr{nu_k}}
    g_k(B(t_{\ell-1}))\,I_3(\de_\ell^{\ot3}).
  \]
Combining the Leibniz rule \eqref{leibnitz} with $D\big(I_q(f^{\ot q})\big)=qI_{q-1}(f^{\ot (q-1)})f$ for any $f\in\HH$, we have
  \begin{multline}\label{eq:computationD3Bis}
  D^3\ang{\la,\VV_n} = \sum_{k=1}^m\la_k\sum_{\ell=1}^{\flr{nu_k}}
    \bigg[g_k'''(B(t_{\ell-1}))I_3(\de_\ell^{\ot3})\ep_{\ell-1}^{\ot3}
  + 9 g_k''(B(t_{\ell-1}))I_2(\de_\ell^{\ot2})
    \ep_{\ell-1}^{\ot2}\tot\de_\ell\\
  + 18 g_k'(B(t_{\ell-1}))I_1(\de_\ell)
    \ep_{\ell-1}\tot\de_\ell^{\ot2}
    + 6 g_k(B(t_{\ell-1}))\de_\ell^{\ot3}\bigg].
  \end{multline}
Combining relations \eqref{eq:computationPartialPhi1}, \eqref{eq:computationPartialPhi}, \eqref{eq:computationD3}, and \eqref{eq:computationD3Bis} we obtain the following expression:
  \begin{multline}\label{decom}
  \frac{\pa}{\pa\la_p} E\left(e^{i\ang{\la,\VV_n}}\xi\right)
    = i E\bigg(e^{i\ang{\la,\VV_n}}\xi\sum_{j=1}^{\flr{nu_p}}
    g_p'''(B(t_{j-1}))\ang{\ep_{j-1},\de_j}_\HH^3\bigg)\\
  - 6\,E\bigg(e^{i\ang{\la,\VV_n}}\xi\sum_{j=1}^{\flr{nu_p}}
    \sum_{k=1}^m\la_k\sum_{\ell=1}^{\flr{nu_k}}
    g_p(B(t_{j-1}))g_k(B(t_{\ell-1}))
    \ang{\de_\ell,\de_j}_\HH^3\bigg)
  + i \sum_{j=1}^{\flr{nu_p}} r_{j,n},
  \end{multline}
with
  \begin{align}\label{eq:rkn}
  r_{j,n} &= i\sum_{k=1}^m\la_k\sum_{\ell=1}^{\flr{nu_k}}
    E\left(g_p(B(t_{j-1}))g_k'''(B(t_{\ell-1}))
    I_3(\de_\ell^{\ot3})e^{i\ang{\la,\VV_n}}\xi\right)
    \ang{\ep_{\ell-1},\de_j}_\HH^3
    \nonumber\displaybreak[0]\\
  & + 9i\sum_{k=1}^m\la_k\sum_{\ell=1}^{\flr{nu_k}}
    E\left(g_p(B(t_{j-1}))g_k''(B(t_{\ell-1}))
    I_2(\de_\ell^{\ot2})e^{i\ang{\la,\VV_n}}\xi\right)
    \ang{\ep_{\ell-1},\de_j}_\HH^2\ang{\de_\ell,\de_j}_\HH
    \nonumber\displaybreak[0]\\
  & + 18i\sum_{k=1}^m\la_k\sum_{\ell=1}^{\flr{nu_k}}
    E\left(g_p(B(t_{j-1}))g_k'(B(t_{\ell-1}))
    I_1(\de_\ell)e^{i\ang{\la,\VV_n}}\xi\right)
    \ang{\de_\ell,\de_j}_\HH^2\ang{\ep_{\ell-1},\de_j}_\HH
    \nonumber\displaybreak[0]\\
  & + 3iE\left(g_p''(B(t_{j-1}))e^{i\ang{\la,\VV_n}}\xi
    \ang{D\ang{\la,\VV_n},\de_j}_\HH\right)
    \ang{\ep_{j-1},\de_j}_\HH^2
    \nonumber\displaybreak[0]\\
  & + 3E\left(g_p''(B(t_{j-1}))e^{i\ang{\la,\VV_n}}
    \ang{D\xi,\de_j}_\HH\right)\ang{\ep_{j-1},\de_j}_\HH^2
    \nonumber\displaybreak[0]\\
  & - 3E\left(g_p'(B(t_{j-1}))e^{i\ang{\la,\VV_n}}\xi
    \ang{D\ang{\la,\VV_n},\de_j}_\HH^2\right)
    \ang{\ep_{j-1},\de_j}_\HH
    \nonumber\displaybreak[0]\\
  & + 6iE\left(g_p'(B(t_{j-1}))e^{i\ang{\la,\VV_n}}
    \ang{D\ang{\la,\VV_n},\de_j}_\HH\ang{D\xi,\de_j}_\HH\right)
    \ang{\ep_{j-1},\de_j}_\HH
    \nonumber\displaybreak[0]\\
  & + 3E\left(g_p'(B(t_{j-1}))e^{i\ang{\la,\VV_n}}
    \ang{D^2\xi,\de_j^{\ot2}}_{\HH^{\ot2}}\right)
    \ang{\ep_{j-1},\de_j}_\HH
    \nonumber\displaybreak[0]\\
  & - iE\left(g_p(B(t_{j-1}))e^{i\ang{\la,\VV_n}}\xi
    \ang{D\ang{\la,\VV_n},\de_j}_\HH^3 \right)
    \nonumber\displaybreak[0]\\
  & - 3E\left(g_p(B(t_{j-1}))e^{i\ang{\la,\VV_n}}
    \ang{D\ang{\la,\VV_n},\de_j}_\HH^2
    \ang{D\xi,\de_j}_\HH\right)
    \nonumber\displaybreak[0]\\
  & + 3iE\left(g_p'(B(t_{j-1}))e^{i\ang{\la,\VV_n}}\xi
    \ang{D^2\ang{\la,\VV_n}),\de_j^{\ot2}}_{\HH^{\ot2}}\right)
    \ang{\ep_{j-1},\de_j}_\HH
    \nonumber\displaybreak[0]\\
  & - 3E\left(g_p(B(t_{j-1}))e^{i\ang{\la,\VV_n}}\xi
    \ang{D\ang{\la,\VV_n},\de_j}_\HH
    \ang{D^2\ang{\la,\VV_n},\de_j^{\ot2}}_{\HH^{\ot2}}\right)
    \nonumber\displaybreak[0]\\
  & + 3iE\left(g_p(B(t_{j-1}))e^{i\ang{\la,\VV_n}}
    \ang{D\xi,\de_j}_\HH
    \ang{D^2\ang{\la,\VV_n},\de_j^{\ot2}}_{\HH^{\ot2}}\right)
    \nonumber\displaybreak[0]\\
  & + 3iE\left(g_p(B(t_{j-1}))e^{i\ang{\la,\VV_n}}
    \ang{D\ang{\la,\VV_n},\de_j}_\HH
    \ang{D^2\xi,\de_j^{\ot2}}_{\HH^{\ot2}}\right)
    \nonumber\displaybreak[0]\\
  & + E\left(g_p(B(t_{j-1}))e^{i\ang{\la,\VV_n}}
    \ang{D^3\xi,\de_j^{\ot3}}_{\HH^{\ot3}}\right)
    \nonumber\displaybreak[0]\\
  & = \sum_{a=1}^{15} R_{j,n}^{(a)}.
  \end{align}
Assume for a moment (see Steps 4 to 8 below) that
  \begin{equation}\label{eq:remain}
  \sum_{j=1}^{\flr{nu_p}} r_{j,n}
    \xrightarrow[n\to\infty]{} 0.
  \end{equation}
By Lemma \ref{lemma:Tec1} $(iv)$ and since $e^{i\ang{\la,\VV_n}}$, $\xi$ and $g_p'''$ are bounded, we have
  \[
  \bigg|E\bigg(e^{i\ang{\la,\VV_n}}\xi
    \sum_{j=1}^{\flr{nu_p}} g_p'''(B(t_{j-1}))
    \ang{\ep_{j-1},\de_j}_\HH^3\bigg)
  - E\bigg(e^{i\ang{\la,\VV_n}}\xi
    \times\frac{(-1)}{8n}\sum_{j=1}^{\flr{nu_p}}
    g_p'''(B(t_{j-1}))\bigg)\bigg| \to 0.
  \]
Moreover, by Lebesgue bounded convergence, we have that
  \[
  \bigg|E\bigg(e^{i\ang{\la,\VV_n}}\xi
    \times\frac{(-1)}{8n}\sum_{j=1}^{\flr{nu_p}}
    g_p'''(B(t_{j-1}))\bigg)
  - E\bigg(e^{i\ang{\la,\VV_n}}\xi
    \times\frac{(-1)}8\int_0^{u_p}
    g_p'''(B(s))\,ds\bigg)\bigg| \to 0.
  \]
Finally, since $(B,\VV_n)\to(B,\VV_\infty)$ in $D_\R[0,\infty)\times\R^m$, we have
  \[
  E\bigg(e^{i\ang{\la,\VV_n}}\xi
    \times\frac{(-1)}8\int_0^{u_p} g_p'''(B(s))\,ds\bigg)
    \to E\bigg(e^{i\ang{\la,\VV_\infty}}\xi
    \times\frac{(-1)}8\int_0^{u_p} g_p'''(B(s))\,ds\bigg).
  \]
Putting these convergences together, we obtain:
  \begin{equation}\label{riemann1}
  E\bigg(e^{i\ang{\la,\VV_n}}\xi
    \sum_{j=1}^{\flr{nu_p}} g_p'''(B(t_{j-1}))
    \ang{\ep_{j-1},\de_j}_\HH^3\bigg)
    \to E\bigg(e^{i\ang{\la,\VV_\infty}}\xi
    \times \frac{(-1)}8\int_0^{u_p} g_p'''(B(s))\,ds\bigg).
  \end{equation}
Similarly, let us show that
  \begin{multline}\label{riemann2}
  6E\bigg(e^{i\ang{\la,\VV_n}}\xi
    \sum_{j=1}^{\flr{nu_p}}\sum_{\ell=1}^{\flr{nu_k}}
    g_p(B(t_{j-1}))g_k(B(t_{\ell-1}))
    \ang{\de_\ell,\de_j}_\HH^3\bigg)\\
  \to \ka^2 E\left(e^{i\ang{\la,\VV_\infty}}\xi
    \times\int_0^{u_p\wedge u_k} g_p(B(s))g_k(B(s))\,ds\right).
  \end{multline}
We have, see \eqref{rho} for the definition of $\rho$:
  \begin{align}
  6\sum_{j=1}^{\flr{nu_p}}\sum_{\ell=1}^{\flr{nu_k}}
    g_p(B(t_{j-1}))&g_k(B(t_{\ell-1}))
    \ang{\de_\ell,\de_j}_\HH^3\notag\\
   &= \frac6n\sum_{j=1}^{\flr{nu_p}}\sum_{\ell=1}^{\flr{nu_k}}
     g_p(B(t_{j-1}))g_k(B(t_{\ell-1}))\rho^3(\ell - j)\notag\\
   &= \frac6n\sum_{r=1-\flr{nu_p}}^{\flr{nu_k}-1}\rho^3(r)
     \sum_{j=1\vee(1-r)}^{\flr{nu_p}\wedge(\flr{nu_k}-r)}
     g_p(B(t_{j-1}))g_k(B(t_{r+j-1})).\label{ascv}
  \end{align}
For each fixed integer $r>0$ (the case $r\le 0$ being similar), we have
  \begin{multline*}
  \bigg|\frac1n\sum_{j=1\vee(1-r)}^{\flr{nu_p}\wedge(\flr{nu_k}-r)}
    g_p(B(t_{j-1}))g_k(B(t_{r+j-1}))
    - \frac1n\sum_{j=1\vee(1-r)}^{\flr{nu_p}\wedge(\flr{nu_k}-r)}
    g_p(B(t_{j-1}))g_k(B(t_{j-1}))\bigg|\\
  \le C\|g_p\|_\infty\sup_{1\le j\le\flr{nu_p}}
    \big|g_k(B(t_{r+j-1})) - g_k(B(t_{j-1}))\big|
    \xrightarrow[n\to\infty]{\text{a.s.}} 0
    \quad\text{by Heine's theorem}.
  \end{multline*}
Hence, for all fixed $r\in\ZZ$,
  \[
  \frac1n\sum_{j=1\vee(1-r)}^{\flr{nu_p}\wedge(\flr{nu_k}-r)}
    g_p(B(t_{j-1}))g_k(B(t_{r+j-1}))
    \xrightarrow[n\to\infty]{\text{a.s.}}
    \int_0^{u_p\wedge u_k}g_p(B(s))g_k(B(s))ds.
  \]
By combining a bounded convergence argument with \eqref{ascv} (observe in particular that $\ka^2=6\sum_{r\in\ZZ}\rho^3(r)<\infty$), we deduce that
  \[
  6\sum_{j=1}^{\flr{nu_p}}\sum_{\ell=1}^{\flr{nu_k}}
    g_p(B(t_{j-1}))g_k(B(t_{\ell-1}))
    \ang{\de_\ell,\de_j}_\HH^3
    \xrightarrow[n\to\infty]{\text{a.s.}}
    \ka^2\int_0^{u_p\wedge u_k}g_p(B(s))g_k(B(s))\,ds.
  \]
Since $(B,\VV_n)\to(B,\VV_\infty)$ in $D_\R[0,\infty)\times\R^m$, we deduce that
  \begin{multline*}
  \bigg(\VV_n,\xi,
    6\sum_{j=1}^{\flr{nu_p}}\sum_{\ell=1}^{\flr{nu_k}}
    g_p(B(t_{j-1}))g_k(B(t_{\ell-1}))
    \ang{\de_\ell,\de_j}_\HH^3\bigg)_{k=1,\ldots,m}\\
  \xrightarrow[]{\text{Law}}
    \bigg(\VV_\infty,\xi,\ka^2\int_0^{u_p\wedge u_k}
    g_p(B(s))g_k(B(s))\,ds\bigg)_{k=1,\ldots,m}
    \quad\text{in $\R^d\times\R\times\R^m$}.
  \end{multline*}
By boundedness of $e^{i\ang{\la,\VV_n}}$, $\xi$ and $g_i$, we have that \eqref{riemann2} follows. Putting \eqref{eq:remain}, \eqref{riemann1}, and \eqref{riemann2} into \eqref{decom}, we deduce \eqref{todo}.

Now, it remains to prove \eqref{eq:remain}.

\bigskip

{\it Step 4.-} Study of $R_{j,n}^{(5)}$, $R_{j,n}^{(8)}$ and $R_{j,n}^{(15)}$ in \eqref{eq:rkn}. Let $k \in \{1,2,3\}$. Since
  \[
  D^k\xi = \sum_{i_1,\ldots,i_k=1}^r
    \frac{\pa^k\psi}{\pa s_{i_1}\cdots\pa s_{i_k}}
    (B_{s_1},\ldots,B_{s_r}){\bf 1}_{[0,s_{i_1}]}\ot
    \ldots\ot{\bf 1}_{[0,s_{i_k}]},
  \]
with $\psi\in C^\infty_b(\R^r)$, we have $\sum_{j=1}^{\flr{nt}}| \ang{D^k\xi,\de_j^{\ot k}}_\HH|\le C n^{-(k-1)/3}$ by Lemma \ref{lemma:Tec1} $(i)$ and $(ii)$. Moreover, $|\ang{\ep_{j-1},\de_j}_\HH|\le n^{-1/3}$ by Lemma \ref{lemma:Tec1} $(i)$. Hence, $\sum_{j=1}^{\flr{nt}}|R_{j,n}^{(p)}| = O(n^{-2/3})\xrightarrow[n\to\infty]{}0$ for $p \in \{5, 8, 15\}$.

\bigskip

{\it Step 5.-} Study of $R_{j,n}^{(2)}$ and $R_{j,n}^{(3)}$ in \eqref{eq:rkn}. We can write, using Lemma \ref{lemma:Tec1}  $(i)$, Cauchy-Schwarz inequality and the definition \eqref{rho} of $\rho$ among other things:
  \begin{multline*}
  \sum_{j=1}^{\flr{nu_p}} \big|R_{j,n}^{(3)}\big|\\
  \le 18\sum_{k=1}^m |\la_k|\sum_{j=1}^{\flr{nu_p}}
    \sum_{\ell=1}^{\flr{nu_k}}
    \big|E\left(g_p(B(t_{j-1}))g'_k(B(t_{\ell-1}))
    I_1(\de_\ell)e^{i\ang{\la,\VV_n}}\xi\right)\big|
    \ang{\de_\ell,\de_j}^2_\HH
    \big|\ang{\ep_{\ell-1},\de_j}_\HH\big|\\
  \le Cn^{-7/6}\sum_{j=1}^{\flr{nu_p}}
    \sum_{\ell=1}^{\flr{nu_k}}\rho(\ell - j)^2
    \le Cn^{-1/6} \sum_{r\in\ZZ}\rho(r)^2 = Cn^{-1/6}
    \xrightarrow[n\to\infty]{} 0.
  \end{multline*}
Concerning $R_{j,n}^{(2)}$, we can write similarly:
  \begin{multline*}
  \sum_{j=1}^{\flr{nu_p}} \big|R_{j,n}^{(2)}\big|\\
  \le 9\sum_{k=1}^m |\la_k|\sum_{j=1}^{\flr{nu_p}}
    \sum_{\ell=1}^{\flr{nu_k}}
    \big|E\left(g_p(B(t_{j-1}))g''_k(B(t_{\ell-1}))
    I_2(\de_\ell^{\ot2})e^{i\ang{\la,\VV_n}}\xi\right)\big|
    \big|\ang{\de_\ell,\de_j}_\HH\big|
    \ang{\ep_{\ell-1},\de_j}_\HH^2\\
  \le Cn^{-4/3}\sum_{j=1}^{\flr{nu_p}}
    \sum_{\ell=1}^{\flr{nu_k}}\big|\rho(\ell - j)\big|
    \le Cn^{-1/3}\sum_{r\in\ZZ}\big|\rho(r)\big| = Cn^{-1/3}
    \xrightarrow[n\to\infty]{} 0.
  \end{multline*}

\bigskip

{\it Step 6.-} Study of $R_{j,n}^{(1)}$, $R_{j,n}^{(6)}$, $R_{j,n}^{(10)}$, and $R_{j,n}^{(12)}$. First, let us deal with $R_{j,n}^{(1)}$. In order to lighten the notation, we set $\wt\xi_{j,\ell}=g_p(B(t_{j-1})) g'''_k(B(t_{\ell-1}))\xi$. Using $I_3(\de_\ell^{\ot3})=I_2(\de_l^{\ot2}) I_1(\de_\ell)-2n^{-1/3}I_1(\de_\ell)$ and then integrating by parts through \eqref{ipp}, we get
  \begin{align*}
  E\left(e^{i\ang{\la,\VV_n}}\wt\xi_{j,\ell}
    I_3(\de_\ell^{\ot3})\right)
    &= E\left(e^{i\ang{\la,\VV_n}}\wt\xi_{j,\ell}
    I_2(\de_\ell^{\ot3})I_1(\de_\ell)\right)
    - 2n^{-1/3} E\left(e^{i\ang{\la,\VV_n}}\wt\xi_{j,\ell}
    I_1(\de_\ell)\right)\\
  &= E\left(e^{i\ang{\la,\VV_n}}
    I_2(\de_\ell^{\ot2})\ang{D\wt\xi_{j,\ell},\de_\ell}_\HH\right)\\
  &\quad + i\la E\left(e^{i\ang{\la,\VV_n}}
    I_2(\de_\ell^{\ot2})\wt\xi_{j,\ell}
    \ang{D\ang{\la,\VV_n},\de_\ell}_\HH\right).
  \end{align*}
Due to Lemma \ref{lemma:Tec1} $(i)$ and Cauchy-Schwarz inequality, we have
  \[
  \sup_{\ell=1,\ldots,\flr{nu_k}}\sup_{j=1,\ldots,\flr{nu_p}}
    \left|E\left(e^{i\ang{\la,\VV_n}}I_2(\de_\ell^{\ot2})
    \ang{D\wt\xi_{j,\ell},\de_\ell}_\HH\right)\right|
    \le Cn^{-2/3}.
  \]
By \eqref{est1} and Cauchy-Schwarz inequality, we also have
  \[
  \sup_{\ell=1,\ldots,\flr{nu_k}}\sup_{j=1,\ldots,\flr{nu_p}}
    \left|E\left(e^{i\ang{\la,\VV_n}}
    I_2(\de_\ell^{\ot2})\wt\xi_{j,\ell}
    \ang{D\ang{\la,\VV_n},\de_\ell}_\HH\right)\right|
    \le Cn^{-2/3}.
  \]
Hence, combined with Lemma \ref{lemma:Tec1} $(i)$ and $(ii)$, we get:
  \begin{multline*}
  \sum_{j=1}^{\flr{nu_p}}\big|R_{j,n}^{(1)}\big|
    \le \sum_{k=1}^m|\la_k|\sum_{j=1}^{\flr{nu_p}}
    \sum_{\ell=1}^{\flr{nu_k}}\big|E\left(
    g_p(B(t_{j-1}))g'''_k(B(t_{\ell-1}))
    I_3(\de_\ell^{\ot3})e^{i\ang{\la,\VV_n}}\xi\right)\big|
    \big|\ang{\ep_{\ell-1},\de_j}_\HH\big|^3\\
  \leq Cn^{-1/3}\sup_{\ell=1,\ldots,\flr{nu_k}}
    \sum_{j=1}^{\flr{nu_p}}
    \big|\ang{\ep_{\ell-1},\de_j}_\HH\big|
    \le Cn^{-1/3} \xrightarrow[n\to\infty]{} 0.
  \end{multline*}
Now, let us concentrate on $R^{(6)}_{j,n}$. Since $e^{i\ang{\la,\VV_n}}$, $\xi$, and $g'_p$ are bounded, we have that
  \begin{align*}
  \sum_{j=1}^{\flr{nu_p}} |R_{j,n}^{(6)}|
    &\le C\sum_{j=1}^{\flr{nu_p}}
    E\left(\ang{D\ang{\la,\VV_n},\de_j}_\HH^2\right)
    |\ang{\ep_{j-1},\de_j}_\HH|\\
  &\le C n^{-1/3}\sum_{j=1}^{\flr{nu_p}}
    E\left(\ang{D\ang{\la,\VV_n},\de_j}_\HH^2\right)
    \quad\text{by Lemma \ref{lemma:Tec1} $(i)$}\\
  &\le C n^{-2/3} \xrightarrow[n\to\infty]{} 0
    \quad\text{by \eqref{est2}}.
  \end{align*}
Similarly,
  \begin{align*}
  \sum_{j=1}^{\flr{nu_p}}|R_{j,n}^{(12)}|
    &\le 3\sum_{j=1}^{\flr{nu_p}} \left\vert E\left(
    g_p(B(t_{j-1}))e^{i\ang{\la,\VV_n}}\xi
    \ang{D\ang{\la,\VV_n},\de_j}_\HH
    \ang{D^2\ang{\la,\VV_n},\de_j^{\ot2}}_{\HH^{\ot2}}
    \right)\right\vert\\
  &\le C\sum_{j=1}^{\flr{nu_p}}
    \bigg(E\left(\ang{D\ang{\la,\VV_n},\de_j}_\HH^2\right)
    + E\left(\ang{D^2\ang{\la,\VV_n},\de_j^{\ot2}}_{\HH^{\ot2}}^2
    \right)\bigg)\\
  &\le C n^{-1/3} \xrightarrow[n\to\infty]{} 0
    \quad\text{by \eqref{est2} and \eqref{est3}}.
  \end{align*}
For $R_{j,n}^{(10)}$, we can write:
  \begin{align*}
  \sum_{j=1}^{\flr{nu_p}}|R_{j,n}^{(10)}|
    &\le 3\sum_{j=1}^{\flr{nu_p}}\left\vert E\left(
    g_p(B(t_{j-1}))e^{i\ang{\la,\VV_n}}
    \ang{D\ang{\la,\VV_n},\de_j}_\HH^2
    \ang{D\xi,\de_j}_\HH\right)\right\vert\\
  &\le Cn^{-1/3}\sum_{j=1}^{\flr{nu_p}}
    \left\vert E\left(
    \ang{D\ang{\la,\VV_n},\de_j}_\HH^2\right)\right\vert
    \quad\text{by Lemma \ref{lemma:Tec1} $(i)$}\\
  &\le Cn^{-2/3} \xrightarrow[n\to\infty]{} 0
    \quad\text{by \eqref{est2}}.
  \end{align*}

\bigskip

{\it Step 7.-} Study of $R^{(4)}_{j,n}$, $R^{(7)}_{j,n}$, $R^{(11)}_{j,n}$, $R^{(13)}_{j,n}$, and $R^{(14)}_{j,n}$.

Using \eqref{eq:DLambdaVn}, and then Cauchy-Schwarz inequality and Lemma \ref{lemma:Tec1} $(i)$, we can write
  \begin{align*}
  &\sum_{j=1}^{\flr{nu_p}}\big|R^{(4)}_{j,n}\big|\\
  &\le 3\sum_{j=1}^{\flr{nu_p}}\left|E\big(
    g''_p(B(t_{j-1})e^{i\ang{\la,\VV_n}}\xi
    \ang{D\ang{\la,\VV_n},\de_j}_\HH\big)\right|
    \ang{\ep_{j-1},\de_j}_\HH^2\\
  &\le 3\sum_{k=1}^m |\la_k|
    \sum_{j=1}^{\flr{nu_p}}\sum_{i=1}^{\flr{nu_k}}\left|E\big(
    g''_p(B(t_{j-1}))g'_k(B(t_{i-1}))e^{i\ang{\la,\VV_n}}\xi
    I_3(\de_i^{\ot3})\big)\right|
    \big|\ang{\ep_{i-1},\de_j}_\HH\big|
    \ang{\ep_{j-1},\de_j}_\HH^2\\
  &\quad + 9\sum_{k=1}^m |\la_k|
     \sum_{j=1}^{\flr{nu_p}}\sum_{i=1}^{\flr{nu_k}}\left|E\big(
     g''_p(B(t_{j-1}))g_k(B(t_{i-1}))e^{i\ang{\la,\VV_n}}\xi
     I_2(\de_i^{\ot2})\big)\right|
     \big|\ang{\de_i,\de_j}_\HH\big|
     \ang{\ep_{j-1},\de_j}_\HH^2\\
  &\le Cn^{-1/6}\sup_{1\le i\le\flr{nu_k}}
     \sum_{j=1}^{\flr{nu_p}}\big|\ang{\ep_{i-1},\de_j}_\HH\big|
     + Cn^{-4/3}\sum_{j=1}^{\flr{nu_p}}
     \sum_{i=1}^{\flr{nu_k}}\big|\rho(i - j)\big|
     \le Cn^{-1/6}\xrightarrow[n\to\infty]{} 0.
 \end{align*}
Using the same arguments, we show that $\sum_{j=1}^{\flr{nu_p}} \big| R^{(7)}_{j,n}\big|\xrightarrow[n\to\infty]{}0$ and $\sum_{j=1}^{\flr{nu_p}} \big|R^{(14)}_{j,n}\big|\xrightarrow[n\to\infty]{}0$.

Differentiating two times in \eqref{eq:DLambdaVn}, we get
  \begin{align*}
  \ang{D^2\ang{\la,\VV_n},\de_j^{\ot2}}_{\HH^{\ot2}}
    &= \sum_{k=1}^m\la_k \sum_{i=1}^{\flr{nu_k}}
    g''_k(B(t_{i-1}))I_3(\de_i^{\ot3})
    \ang{\ep_{i-1},\de_j}^2_\HH\\
  &\quad + 6\sum_{k=1}^m\la_k \sum_{i=1}^{\flr{nu_k}}
    g'_k(B(t_{i-1}))I_2(\de_i^{\ot2})
    \ang{\ep_{i-1},\de_j}_\HH\ang{\de_i,\de_j}_\HH\\
  &\quad + 6\sum_{k=1}^m\la_k \sum_{i=1}^{\flr{nu_k}}
    g_k(B(t_{i-1}))I_1(\de_i)\ang{\de_i,\de_j}_\HH^2.
  \end{align*}
Hence, using Cauchy-Schwarz inequality and Lemma \ref{lemma:Tec1} $(i)$-$(ii)$, we can write
  \begin{align*}
  &\sum_{j=1}^{\flr{nu_p}}\big|R^{(11)}_{j,n}\big|\\
  &\le 3\sum_{j=1}^{\flr{nu_p}}\left|E\big(
    g'_p(B(t_{j-1}))e^{i\ang{\la,\VV_n}}\xi
    \ang{D^2\ang{\la,\VV_n},\de_j^{\ot2}}_{\HH^{\ot2}}\big)\right|
    \ang{\ep_{j-1},\de_j}_\HH\displaybreak[0]\\
  &\le 3\sum_{k=1}^m |\la_k|
    \sum_{j=1}^{\flr{nu_p}}\sum_{i=1}^{\flr{nu_k}}\left|E\big(
    g'_p(B(t_{j-1}))g''_k(B(t_{i-1}))e^{i\ang{\la,\VV_n}}\xi
    I_3(\de_i^{\ot3})\big)\right|
    \ang{\ep_{i-1},\de_j}_\HH^2
    \big|\ang{\ep_{j-1},\de_j}_\HH\big|\displaybreak[0]\\
  & + 18\sum_{k=1}^m |\la_k|
    \sum_{j=1}^{\flr{nu_p}}\sum_{i=1}^{\flr{nu_k}}\left|E\big(
    g'_p(B(t_{j-1}))g'_k(B(t_{i-1}))e^{i\ang{\la,\VV_n}}\xi
    I_2(\de_i^{\ot2})\big)\right|\\
  & \hskip9cm\times\big|\ang{\ep_{i-1},\de_j}_\HH\big|
    \big|\ang{\ep_{j-1},\de_j}_\HH\big|
    \big|\ang{\de_i,\de_j}_\HH\big|\displaybreak[0]\\
  & + 18\sum_{k=1}^m |\la_k|
    \sum_{j=1}^{\flr{nu_p}}\sum_{i=1}^{\flr{nu_k}}\left|E\big(
    g'_p(B(t_{j-1}))g_k(B(t_{i-1}))e^{i\ang{\la,\VV_n}}\xi
    I_1(\de_i)\big)\right|
    \big|\ang{\ep_{j-1},\de_j}_\HH\big|
    \ang{\de_i,\de_j}_\HH^2\displaybreak[0]\\
  &\le Cn^{-1/6}\sup_{1\le i\le\flr{nu_k}}
    \sum_{j=1}^{\flr{nu_p}}\big|\ang{\ep_{i-1},\de_j}_\HH\big|
    + Cn^{-4/3}\sum_{j=1}^{\flr{nu_p}}\sum_{i=1}^{\flr{nu_k}}
    \big|\rho(i - j)\big|\\
  &\hskip9cm + Cn^{-5/6}\sum_{j=1}^{\flr{nu_p}}
    \sum_{i=1}^{\flr{nu_k}} \big|\rho(i - j)\big|\displaybreak[0]\\
  &\le Cn^{-1/6} \xrightarrow[n\to\infty]{} 0.
  \end{align*}
Using the same arguments, we show that $\sum_{j=1}^{\flr{nu_p}}\big| R^{(13)}_{j,n}\big|\xrightarrow[n \to \infty]{}0$.

\bigskip

{\it Step 8.-}
Now, we consider the last term in \eqref{eq:rkn}, that is $R_{j,n}^{(9)}$. Since $e^{i\ang{\la,\VV_n}}$, $\xi$, and $g_p$ are bounded, we can write
  \begin{align*}
  \bigg|\sum_{j=1}^{\flr{nu_p}} R_{j,n}^{(9)}\bigg|
    &\le C\sum_{j=1}^{\flr{nu_p}} E\big(
    |\ang{D\ang{\la,\VV_n},\de_j}_\HH|^3\big)\\
  &\le C\sum_{j=1}^{\flr{nu_p}} E\big(
    \ang{D\ang{\la,\VV_n},\de_j}_\HH^2\big)
    + E\big(\ang{D\ang{\la,\VV_n},\de_j}_\HH^4\big).
  \end{align*}
In addition we have, see \eqref{eq:DLambdaVn}, that
  \begin{align*}
  \sum_{j=1}^{\flr{nu_p}} &E(\ang{D\ang{\la,\VV_n},\de_j}_\HH^4)\\
  &\le 8m^3 \sum_{j=1}^{\flr{nu_p}}\sum_{k=1}^m
    \la_k^4\sum_{i_1=1}^{\flr{nu_k}}\sum_{i_2=1}^{\flr{nu_k}}
    \sum_{i_3=1}^{\flr{nu_k}}\sum_{i_4=1}^{\flr{nu_k}}
    E\bigg(\prod_{a=1}^4\ang{\ep_{t_{i_a}},\de_j}_\HH
    g_k'(B(t_{i_1-1}))I_3(\de_{i_a}^{\ot3})\bigg)\\
  &\quad  + 648m^3\sum_{j=1}^{\flr{nu_p}}\sum_{k=1}^m\la_k^4
    \sum_{i_1=1}^{\flr{nu_k}}\sum_{i_2=1}^{\flr{nu_k}}
    \sum_{i_3=1}^{\flr{nu_k}}\sum_{i_4=1}^{\flr{nu_k}}
    E\bigg(\prod_{a=1}^4\ang{\de_{i_a},\de_j}_\HH
    g_k(B(t_{i_a-1}))I_2(\de_{i_a}^{\ot2})\bigg)\\
  &\le C\sum_{k=1}^m\la_k^4\bigg[n^{-4/3}\sum_{j=1}^{\flr{nu_p}}
    \sum_{i_1=1}^{\flr{nu_k}}\sum_{i_2=1}^{\flr{nu_k}}
    \sum_{i_3=1}^{\flr{nu_k}}\sum_{i_4=1}^{\flr{nu_k}}
    \bigg|E\bigg(\prod_{a=1}^4
    g_k'(B(t_{i_1-1}))I_3(\de_{i_a}^{\ot3})\bigg)\bigg|\\
  &\quad + \sum_{j=1}^{\flr{nu_p}}
    \sum_{i_1=1}^{\flr{nu_k}}\sum_{i_2=1}^{\flr{nu_k}}
    \sum_{i_3=1}^{\flr{nu_k}}\sum_{i_4=1}^{\flr{nu_k}}
    \prod_{a=1}^4 |\ang{\de_{i_a},\de_j}_\HH|
    \bigg|E\bigg(\prod_{a=1}^4
    g_k(B(t_{i_a-1}))I_2(\de_{i_a}^{\ot2})\bigg)\bigg|\bigg].
  \end{align*}
By Lemma \ref{lemma:Tec3} we have that
  \[
  \sum_{i_1=1}^{\flr{nu_k}}\sum_{i_2=1}^{\flr{nu_k}}
    \sum_{i_3=1}^{\flr{nu_k}}\sum_{i_4=1}^{\flr{nu_k}}
    \bigg|E\bigg(\prod_{a=1}^4
    g_k'(B(t_{i_1-1}))I_3(\de_{i_a}^{\ot3})\bigg)\bigg|
    \le C,
  \]
so that
  \[
  n^{-4/3}\sum_{j=1}^{\flr{nu_p}}
    \sum_{i_1=1}^{\flr{nu_k}}\sum_{i_2=1}^{\flr{nu_k}}
    \sum_{i_3=1}^{\flr{nu_k}}\sum_{i_4=1}^{\flr{nu_k}}
    \bigg|E\bigg(\prod_{a=1}^4
    g_k'(B(t_{i_1-1}))I_3(\de_{i_a}^{\ot3})\bigg)\bigg|
    \le C n^{-1/3}.
  \]
On the other hand, by Cauchy-Schwarz inequality, we have
  \[
  \bigg|E\bigg(\prod_{a=1}^4
    g_k(B(t_{i_a-1}))I_2(\de_{i_a}^{\ot2})\bigg)\bigg|
    \le C,
  \]
so that, with $\rho$ defined by \eqref{rho},
  \begin{multline*}
  \sum_{j=1}^{\flr{nu_p}}
    \sum_{i_1=1}^{\flr{nu_k}}\sum_{i_2=1}^{\flr{nu_k}}
    \sum_{i_3=1}^{\flr{nu_k}}\sum_{i_4=1}^{\flr{nu_k}}
    \prod_{a=1}^4 |\ang{\de_{i_a},\de_j}_\HH|
    \bigg|E\bigg(\prod_{a=1}^4 g_k(B(t_{i_a-1}))
    I_2(\de_{i_a}^{\ot2})\bigg)\bigg|\\
  \le Cn^{-4/3}\sum_{j=1}^{\flr{nu_p}}
    \sum_{i_1=1}^{\flr{nu_k}} |\rho(i_1 - j)|
    \times \sum_{i_2=1}^{\flr{nu_k}} |\rho(i_2 - j)|
    \times \sum_{i_3=1}^{\flr{nu_k}} |\rho(i_3 - j)|
    \times \sum_{i_4=1}^{\flr{nu_k}} |\rho(i_4 - j)|\\
    \le Cn^{-1/3}\bigg(\sum_{r\in\ZZ}|\rho(r)|\bigg)^4
    = Cn^{-1/3}.
\end{multline*}
As a consequence, combining the previous estimates with \eqref{est2}, we have shown that
  \[
  \bigg|\sum_{j=1}^{\flr{nu_p}} R_{j,n}^{(9)}\bigg|
    \le Cn^{-1/3} \xrightarrow[n\to\infty]{} 0,
  \]
and the proof of Theorem \ref{fdd2} is done. \qed

\begin{thm}\label{T:fdd3}
If $g\in C^\infty(\RR)$ is bounded with bounded derivatives, then the sequence $(B,G_n^+(1,B),G_n^+(g,B))$ converges to $(B, \cub{B},(1/8)\int g'''(B)\,ds+\int g(B)d\cub{B})$ in the sense of finite-dimensional distributions on $[0,\infty)$, where
  \[
  G_n^+(g,B,t) := \frac1{\sqrt{n}}\sum_{j=1}^{\flr{nt}}
    g(B(t_j))\,h_3(n^{1/6}\De B_j), \quad t\ge 0, \quad n\ge 1.
  \]
\end{thm}

\pf The proof is exactly the same as the proof of Theorem \ref{fdd2}, except $\ep_{j-1}$ must be everywhere replaced $\ep_j$, and Lemma \ref{lemma:Tec1} $(v)$ must used instead of Lemma \ref{lemma:Tec1} $(iv)$. \qed

\bigskip

\noindent{\bf Proof of Theorem \ref{T:main_fdd}.} We begin by observing the following general fact. Suppose $U$ and $V$ are c\`adl\`ag processes adapted to a filtration under which $V$ is a semimartingale. Similarly, suppose $\wt U$ and $\wt V$ are c\`adl\`ag processes adapted to a filtration under which $\wt V$ is a semimartingale. If the processes $(U,V)$ and $(\wt U,\wt V)$ have the same law, then $\int _0^\cdot U(s-)\,dV(s)$ and $\int_0^\cdot\wt U(s-)\,d\wt V(s)$ have the same law. This is easily seen by observing that these integrals are the limit in probability of left-endpoint Riemann sums.

Now, let $G_n^-$ and $G_n^+$ be as defined previously in this section. Define
  \begin{align*}
  G^-(g,B,t) &= -\frac18\int_0^t g'''(B(s))\,ds
    + \int_0^t g(B(s))\,d\cub{B}_s,\\
  G^+(g,B,t) &= \frac18\int_0^t g'''(B(s))\,ds
    + \int_0^t g(B(s))\,d\cub{B}_s.
  \end{align*}
Let ${\bf t}=(t_1,\ldots, t_d)$, where $0\le t_1<\cdots<t_d$. Let
  \[
  G_n^-(g,B,{\bf t}) = (G_n^-(g,B,t_1), \ldots, G_n^-(g,B,t_d)),
  \]
and similarly for $G_n^+$, $G^-$, and $G^+$. By Theorems \ref{T:signed_cubic}, \ref{fdd2}, and \ref{T:fdd3}, the sequence
  \[
  \{(B, V_n(B), G_n^-(g,B,{\bf t}), G_n^+(g,B,{\bf t}))\}_{n=1}^\infty
  \]
is relatively compact in $D_{\RR^2}[0,\infty)\times\RR^d\times\RR^d$. By passing to a subsequence, we may assume it converges in law in $D_{\RR^2} [0,\infty)\times\RR^d\times\RR^d$ to $(B,\cub{B},X,Y)$, where $X,Y\in\RR^d$.

By Theorems \ref{T:signed_cubic} and \ref{fdd2}, $\{(B, V_n(B), G_n^-(g, B,{\bf t}))\}$ is relatively compact in $D_{\RR^2}[0,\infty)\times\RR^d$, and converges in the sense of finite-dimensional distributions to $(B,\cub{B},G^-(g,B,{\bf t}))$. It follows that $(B, V_n(B), G_n^-(g, B,{\bf t}))\to(B,\cub{B},G^-(g,B,{\bf t}))$ in law in $D_{\RR^2}[0, \infty)\times\RR^d$. Hence, $(B,\cub{B},X)$ and $(B,\cub{B},G^-(g, B,{\bf t}))$ have the same law in $D_{\RR^2}[0,\infty)\times\RR^d$. By the general fact we observed at the beginning of the proof, $(G^-(g, B),X)$ and $(G^-(g,B),G^-(g,B,{\bf t}))$ have the same law. In particular, $(G^-(g,B,{\bf t}),X)$ and $(G^-(g,B,{\bf t}),G^-(g,B,{\bf t}))$ have the same law. But this implies $G^-(g,B,{\bf t})-X$ has the same law as the zero random variable, which gives $G^-(g,B,{\bf t})-X=0$ a.s.

We have thus shown that $X=G^-(g,B,{\bf t})$ a.s. Similarly, $Y = G^+(g,B,{\bf t})$ a.s. It follows that
  \[
  (B, V_n(B), G_n^-(g,B), G_n^+(g,B))
    \to (B, \cub{B}, G^-(g,B), G^+(g,B)),
  \]
and therefore
  \[
  \left(B,V_n(B),\frac{G_n^-(g,B) + G_n^+(g,B)}2\right)
    \to \left(B,[\![B]\!],\frac{G^-(g,B) + G^+(g,B)}2\right),
  \]
in the sense of finite-dimensional distributions on $[0,\infty)$, which is what was to be proved. \qed

\section{Moment bounds}

The following four moment bounds are central to our proof of relative compactness in Theorem \ref{T:main}.

\begin{thm}\label{T:4th_mom}
There exists a constant $C$ such that
  \[
  E|V_n(B,t) - V_n(B,s)|^4 \le C\left({\frac{\flr{nt} - \flr{ns}}{n}
    }\right)^2,
  \]
for all $n$, $s$, and $t$.
\end{thm}

\pf The calculations in the proof of Theorem 10 in \cite{NO} show that
  \[
  E\bigg|\sum_{j=\flr{ns}+1}^{\flr{nt}}\De B_j^3\bigg|^{2p}
    \le C_p\left({\frac{\flr{nt} - \flr{ns}}{n}}\right)^p,
  \]
for all $n$, $s$, and $t$. \qed

\begin{thm}\label{T:5int_mom}
Let $g\in C^1(\RR)$ have compact support. Fix $T>0$ and let $c$ and $d$ be integers such that $0\le t_c<t_d\le T$. Then
  \[
  E\bigg|\sum_{j=c+1}^d g(\be_j)\De B_j^5\bigg|^2
    \le C\|g\|_{1,\infty}^2\De t^{1/3}|t_d - t_c|^{4/3},
  \]
where $\|g\|_{1,\infty} = \|g\|_\infty+\|g'\|_\infty$, and $C$ depends only on $T$.
\end{thm}

\pf Note that
  \begin{equation}\label{5int_mom.1}
  E\bigg|\sum_{j=c+1}^d g(\be_j)\De B_j^5\bigg|^2
    = \sum_{i=c+1}^d\sum_{j=c+1}^d E_{ij},
  \end{equation}
where $E_{ij} = E[g(\be_i)\De B_i^5 g(\be_j)\De B_j^5]$. Let $K = \|g\|_{1,\infty}$, and define $f:\RR^3\to\RR$ by $f(x)=K^{-2}g(x_1) g(x_2)x_3^5$. Note that $f$ has polynomial growth of order 1 with constants $K=1$ and $r=5$.

Let $\xi_1=\be_i$, $\xi_2=\be_j$, $\xi_3=\De t^{-1/6}\De B_i$, $Y =\De t^{-1/6} \De B_j$, and $\ph(y)=y^5$. Then $E_{ij}=K^2\De t^{5/3} E[f(\xi)\ph(Y)]$. By Theorem \ref{T:Gauss_Taylor} with $k=0$, $|E[f(\xi) \ph(Y)]| \le C|\eta|$, where $\eta_j =E[\xi_jY]$. Using Lemma \ref{L:covar}, we have
  \begin{align*}
  |\eta_1| &\le C\De t^{1/6}(j^{-2/3} + |j - i|_+^{-2/3}),\\
  |\eta_2| &\le C\De t^{1/6}j^{-2/3},\\
  |\eta_3| &\le C|j- i|_+^{-5/3}.
  \end{align*}
Hence,
  \[
  |E_{ij}| = K^2\De t^{5/3}|E[f(\xi)\ph(Y)]|
    \le C K^2(\De t^{11/6}(j^{-2/3} + |j - i|_+^{-2/3})
    + \De t^{5/3}|j - i|_+^{-5/3}).
  \]
Substituting this into \eqref{5int_mom.1} gives
  \[
  \begin{split}
  E\bigg|\sum_{j=c+1}^d g(B(t_j))\De B_j^5\bigg|^2
    &\le C K^2(\De t^{11/6}(d - c)^{4/3} + \De t^{5/3}(d - c))\\
  &\le C K^2\De t^{5/3}(d - c)^{4/3}
    = C K^2\De t^{1/3}|t_d - t_c|^{4/3},
  \end{split}
  \]
which completes the proof. \qed

\begin{thm}\label{T:3int_mom.N}
Let $g\in C^2(\RR)$ have compact support. Fix $T>0$ and let $c$ and $d$ be integers such that $0\le t_c<t_d\le T$. Then
  \[
  E\bigg|\sum_{j=c+1}^d g(\be_j)\De B_j^3\bigg|^2
    \le \|g\|_{2,\infty}^2|t_d - t_c|,
  \]
where $\|g\|_{2,\infty} = \|g\|_\infty+\|g'\|_\infty+\|g''\|_\infty$, and $C$ depends only on $T$.
\end{thm}

\pf Note that
  \begin{equation}\label{3int_mom.N.1}
  E\bigg|\sum_{j=c+1}^d g(\be_j)\De B_j^3\bigg|^2
    = \sum_{i=c+1}^d\sum_{j=c+1}^d E_{ij},
  \end{equation}
where $E_{ij} = E[g(\be_i)\De B_i^3 g(\be_j)\De B_j^3]$. Let $K = \|g\|_{2,\infty}$, and define $f: \RR^3\to\RR$ by $f(x) = K^{-2} g(x_1)g(x_2) x_3^3$. Note that $f$ has polynomial growth of order 2 with constants $K=1$ and $r=3$.

Let $\xi_1=\be_i$, $\xi_2=\be_j$, $\xi_3=\De t^{-1/6}\De B_i$, $Y =\De t^{-1/6} \De B_j$, and $\ph(y)=y^3$. Then $E_{ij} =K^2\De t E[f(\xi)\ph(Y)]$. By Theorem \ref{T:Gauss_Taylor} with $k=1$, $E[f(\xi) \ph(Y)]| = \eta_1 E[\pa_1 f(\xi)] + \eta_2 E[\pa_2 f(\xi)] + R$, where $|R|\le C(|\eta_3|+|\eta|^2)$. By \eqref{Gauss_Taylor}, if $j=1$ or $j=2$, $|E[\pa_j f(\xi)]|\le C(|E[\xi_1\xi_3]|+|E[\xi_2\xi_3]|)$. Therefore, using $|\eta_3|^2\le|\eta_3|$ and $|ab|\le|a|^2 +|b|^2$,
  \[
  |E_{ij}| \le CK^2\De t(|\eta_3| + |\eta_1|^2 + |\eta_2|^2
    + |E[\xi_1\xi_3]|^2 + |E[\xi_2\xi_3]|^2).
  \]
Using Lemma \ref{L:covar}, we have
  \begin{align*}
  |E[\xi_2\xi_3]| &\le C\De t^{1/6}(i^{-2/3} + |j - i|_+^{-2/3}),\\
  |E[\xi_1\xi_3]| &\le C\De t^{1/6}i^{-2/3}.
  \end{align*}
Together with the estimates from the proof of Theorem \ref{T:5int_mom}, this gives
  \[
  |E_{ij}| \le CK^2(\De t^{4/3}(i^{-4/3} + j^{-4/3} + |j - i|_+^{-4/3})
    + \De t|j - i|_+^{-5/3}).
  \]
Substituting this into \eqref{3int_mom.N.1} gives
  \[
  E\bigg|\sum_{j=c+1}^d g(B(t_j))\De B_j^3\bigg|^2
    \le CK^2\De t(d - c) = CK^2|t_d - t_c|,
  \]
which completes the proof. \qed

\begin{thm}\label{T:3int_mom}
Suppose $g\in C^3(\RR)$ has compact support. Fix $T>0$ and let $c$ and $d$ be integers such that $0\le t_c<t_d\le T$. Then
  \[
  E\bigg|\sum_{j=c+1}^d (g(\be_j) - g(\be_c))\De B_j^3\bigg|^2
    \le C|t_d - t_c|^{4/3},
  \]
where $C$ depends only on $g$ and $T$.
\end{thm}

\pf Let $Y_j=g(\be_j)-g(\be_c)$, and note that
  \begin{equation}\label{3int_mom.3}
  E\bigg|\sum_{j=c+1}^d Y_j\De B_j^3\bigg|^2
    = \sum_{i=c+1}^d\sum_{j=c+1}^d E_{ij},
  \end{equation}
where $E_{ij}=E[Y_i\De B_i^3 Y_j\De B_j^3]$. For fixed $i,j$, define $f: \RR^4\to\RR$ by
  \[
  f(x) = \left({\frac{g(x_1 + \si_ix_2) - g(x_1)}{\si_i}}\right)
    \left({\frac{g(x_1 + \si_jx_3) - g(x_1)}{\si_j}}\right)x_4^3,
  \]
where $\si_j^2=E|\be_j-\be_c|^2$. Note that $f$ has polynomial growth of order 2 with constants $K$ and $r$ that do not depend on $i$ or $j$.

Let $\xi_1=\be_c$, $\xi_2=\si_i^{-1}(\be_i-\be_c)$, $\xi_3=\si_j^{-1}(\be_j -\be_c)$, $\xi_4=\De t^{-1/6}\De B_i$, $Y=\De t^{-1/6}\De B_j$, and $\ph(y)=y^3$. Note that $E_{ij} = \si_i\si_j\De t E[f(\xi)\ph(Y)]$, so that by Lemma \ref{L:covar}(v),
  \begin{equation}\label{3int_mom.4}
  |E_{ij}| \le C\De t^{4/3}|i - c|^{1/6}|j - c|^{1/6}|E[f(\xi)\ph(Y)]|.
  \end{equation}
By Theorem \ref{T:Gauss_Taylor} with $k=1$, $E[f(\xi)\ph(Y)]=3\sum_{k=1}^3 \eta_k E[\pa_kf(\xi)] + R$, where $|R| \le C(|\eta_4| + |\eta|^2)$. Using $|ab|\le|a|^2 +|b|^2$ and the fact that $|\eta_j|^2\le|\eta|^2$, this gives
  \[
  |E[f(\xi)\ph(Y)]| \le C\bigg(\sum_{k=1}^3|E[\pa_k f(\xi)]|^2
    + |\eta_4| + |\eta|^2\bigg).
  \]
By \eqref{Gauss_Taylor}, for each $k\le3$, $|E[\pa_k f(\xi)]| \le C\sum_{\ell=1}^3 |E[\xi_\ell\xi_4]|$. Therefore, since $\eta_j =E[\xi_jY]$, we have
  \begin{equation}\label{3int_mom.5}
  |E[f(\xi)\ph(Y)]| \le C\bigg(E[\xi_4Y] + \sum_{k=1}^3
    (|E[\xi_kY]|^2 +|E[\xi_k\xi_4]|^2)\bigg).
  \end{equation}
To estimate these covariances, first note that $d-c=n(t_d-t_c)\le nT$. Hence, $\De t=n^{-1}\le C(d-c)^{-1}$. Now, using Lemma \ref{L:covar},
  \begin{align*}
  |E[\xi_1Y]| &\le C\De t^{1/6}|j - c|^{-2/3}
    \le C|d - c|^{-1/6}|j - c|^{-2/3} \le C|j - c|^{-5/6},\\
  |E[\xi_2Y]| &\le C|i - c|^{-1/6}(|j - c|^{-2/3} + |j - i|_+^{-2/3}),\\
  |E[\xi_3Y]| &\le C|j - c|^{-1/6}|j - c|^{-2/3} = C|j - c|^{-5/6},\\
  |E[\xi_4Y]| &\le C|j - i|_+^{-5/3}.
  \end{align*}
Similarly,
  \begin{align*}
  |E[\xi_1\xi_4]| &\le C\De t^{1/6}|i - c|^{-2/3}
    \le C|d - c|^{-1/6}|i - c|^{-2/3} \le C|i - c|^{-5/6},\\
  |E[\xi_2\xi_4]| &\le C|i - c|^{-1/6}|i - c|^{-2/3} = C|i - c|^{-5/6},\\
  |E[\xi_3\xi_4]| &\le C|j - c|^{-1/6}(|i - c|^{-2/3} + |j - i|_+^{-2/3}).
  \end{align*}
Substituting these estimates into \eqref{3int_mom.5} and using \eqref{3int_mom.4} gives
  \begin{multline*}
  |E_{ij}| \le C\De t^{4/3}(|i - c|^{1/6}|j - c|^{1/6}|j - i|_+^{-5/3}\\
  + |i - c|^{1/6}|j - c|^{-3/2} + |i - c|^{-1/6}|j - c|^{-7/6}
    + |i - c|^{-1/6}|j - c|^{1/6}|j - i|_+^{-4/3}\\
  + |i - c|^{-3/2}|j - c|^{1/6} + |i - c|^{-7/6}|j - c|^{-1/6}
    + |i - c|^{1/6}|j - c|^{-1/6}|j - i|_+^{-4/3}).
  \end{multline*}
We can simplify this to
  \begin{multline*}
  |E_{ij}| \le C\De t^{4/3}(|i - c|^{1/6}|j - c|^{1/6}|j - i|_+^{-4/3}\\
  + |i - c|^{1/6}|j - c|^{-7/6}
    + |j - c|^{1/6}|j - i|_+^{-4/3}\\
  + |i - c|^{-7/6}|j - c|^{1/6}
    + |i - c|^{1/6}|j - i|_+^{-4/3}).
  \end{multline*}
Using $|ab|\le|a|^2+|b|^2$, this further simplifies to
  \begin{multline*}
  |E_{ij}| \le C\De t^{4/3}(|i - c|^{1/3}|j - i|_+^{-4/3}
    + |j - c|^{1/3}|j - i|_+^{-4/3}\\
  + |i - c|^{1/6}|j - c|^{-7/6}
  + |i - c|^{-7/6}|j - c|^{1/6}).
  \end{multline*}
We must now make use of \eqref{3int_mom.3}. Note that
  \begin{multline*}
  \De t^{4/3}\sum_{i=c+1}^d\sum_{j=c+1}^d
    |i - c|^{1/3}|j - i|_+^{-4/3}
    \le C\De t^{4/3}\sum_{i=c+1}^d |i - c|^{1/3}\\
    \le C\De t^{4/3}(d - c)^{4/3} = C|t_d - t_c|^{4/3}.
  \end{multline*}
Similarly,
  \[
  \De t^{4/3}\sum_{j=c+1}^d\sum_{i=c+1}^d
    |j - c|^{1/3}|j - i|_+^{-4/3}
    \le C|t_d - t_c|^{4/3}.
  \]
Also,
  \begin{multline*}
  \De t^{4/3}\sum_{i=c+1}^d\sum_{j=c+1}^d
    |i - c|^{1/6}|j - c|^{-7/6}
    \le C\De t^{4/3}\sum_{i=c+1}^d |i - c|^{1/6}
    \le C\De t^{4/3}(d - c)^{7/6}\\
    \le C\De t^{4/3}(d - c)^{4/3} = C|t_d - t_c|^{4/3},
  \end{multline*}
and similarly,
  \[
  \De t^{4/3}\sum_{j=c+1}^d\sum_{i=c+1}^d
    |i - c|^{-7/6}|j - c|^{1/6}
    \le C|t_d - t_c|^{4/3}.
  \]
It follows, therefore, that $\sum_{i=c+1}^d\sum_{j=c+1}^d|E_{ij}|\le C |t_d -t_c|^{4/3}$. By \eqref{3int_mom.3}, this completes the proof. \qed

\section{Proof of main result}\label{S:main_proof}

\begin{lemma}\label{L:5int_conv}
If $g\in C^1(\RR)$ has compact support, then $\sum_{j=1}^{\flr{nt}} g(\be_j)\De B_j^5 \to 0$ ucp.
\end{lemma}

\pf Let $X_n(g,t) = \sum_{j=1}^{\flr{nt}}g(\be_j)\De B_j^5$. Fix $T>0$ and let $0\le s<t\le T$ be arbitrary. Then
  \[
  X_n(g,t) - X_n(g,s) = \sum_{j=c+1}^d g(\be_j)\De B_j^5,
  \]
where $c=\flr{ns}$ and $d=\flr{nt}$. By Theorem \ref{T:5int_mom},
  \[
  E|X_n(g,t) - X_n(g,s)|^2 \le C\De t^{1/3}|t_d - t_c|^{4/3}
    \le C|t_d - t_c|^{5/3}
    = C\left({\frac{\flr{nt} - \flr{ns}}n}\right)^{5/3},
  \]
where $C$ depends only on $g$ and $T$. This verifies condition \eqref{moment_cond} of Theorem \ref{T:moment_cond}. By Theorem \ref{T:5int_mom}, $\sup_n E|X_n(g,T)|^2 \le CT^{4/3}<\infty$. Hence, by Theorem \ref{T:moment_cond}, $\{X_n(g)\}$ is relatively compact in $D_\RR[0,\infty)$. By Lemma \ref{L:in_prob}, it will therefore suffice to show that $X_n(g,t)\to0$ in probability for each fixed $t$. But this follows easily by taking $s=0$ above, which gives $E|X_n (g,t)|^2\le C\De t^{1/3}$ and completes the proof. \qed

\begin{lemma}\label{L:Taylor_expan}
If $g\in C^6(\RR)$ has compact support, then
  \[
  I_n(g',B,t) \approx g(B(t)) - g(B(0))
    + \frac1{12}\sum_{j=1}^{\flr{nt}}g'''(\be_j)\De B_j^3.
  \]
\end{lemma}

\pf Fix $a,b \in\RR$. Let $x=(a+b)/2$ and $h=(b-a)/2$. By Theorem \ref{T:Taylor},
  \[
  \begin{split}
  g(b) - g(a) &= (g(x+h) - g(x)) - (g(x-h) - g(x))\\
  &= \sum_{j=1}^6 g^{(j)}(x)\frac{h^j}{j!}
    - \sum_{j=1}^6 g^{(j)}(x)\frac{(-h)^j}{j!}
    + R_1(x,h) - R_1(x,-h)\\
  &= \sum_{\substack{j=1\\\text{$j$ odd}}}^6 \frac1{j!2^{j-1}}
    g^{(j)}(x)(b - a)^j + R_1(x,h) - R_1(x,-h)\\
  &= g'(x)(b - a) + \frac1{24}g'''(x)(b - a)^3
    + \frac1{5!2^4}g^{(5)}(x)(b - a)^5 + R_2(a,b),
  \end{split}
  \]
where $R_2(a,b)=R_1(x,h)-R_1(x,-h)$ and
  \[
  R_1(x,h) = \frac{h^6}{5!}\int_0^1 (1 - u)^6
    [g^{(6)}(x + uh) - g^{(6)}(x)]\,du.
  \]
Similarly,
  \[
  \begin{split}
  \frac{g'(a) + g'(b)}2 - g'(x) &= \frac12(g'(x+h) - g'(x))
    + \frac12(g'(x-h) - g'(x))\\
  &= \frac12\sum_{j=1}^5 g^{(j+1)}(x)\frac{h^j}{j!}
    + \frac12\sum_{j=1}^5 g^{(j+1)}(x)\frac{(-h)^j}{j!} + R_4(a,b)\\
  &= \frac18 g'''(x)(b - a)^2 + \frac1{4!2^4}g^{(5)}(x)(b - a)^4
    + R_4(a,b),
  \end{split}
  \]
where $R_4(a,b)=R_3(x,h)+R_3(x,-h)$ and
  \[
  R_3(x,h) = \frac{h^5}{4!}\int_0^1 (1 - u)^5
    [g^{(6)}(x + uh) - g^{(6)}(x)]\,du.
  \]
Combining these two expansions gives
  \[
  g(b) - g(a) = \frac{g'(a) + g'(b)}2(b - a)
    - \frac1{12}g'''(x)(b - a)^3
    + \ga g^{(5)}(x)(b - a)^5 + R_6(a,b),
  \]
where $\ga=(5!2^4)^{-1}-(4!2^4)^{-1}$ and
  \[
  R_6(a,b) = R_2(a,b) - R_4(a,b)(b - a).
  \]
Note that $R_6(a,b)=h(a,b)(b - a)^6$, where
  \[
  |h(a,b)| \le C\sup_{0\le u\le 1}|g^{(6)}(x+uh)-g^{(6)}(x)|.
  \]
Taking $a=B(t_{j-1})$ and $b=B(t_j)$ gives
  \begin{multline*}
  g(B(t_j)) - g(B(t_{j-1})) = \frac{g'(B(t_{j-1}))
    + g'(B(t_j))}2\De B_j - \frac1{12}g'''(\be_j)\De B_j^3
    + \ga g^{(5)}(\be_j)\De B_j^5\\
  + h(B(t_{j-1}),B(t_j))\De B_j^6
  \end{multline*}
Recall that $B_n(t)=B(\flr{nt}/n)$, so that
  \[
  g(B(t)) - g(B(0)) = I_n(g',B,t)
    - \frac1{12}\sum_{j=1}^{\flr{nt}} g'''(\be_j)\De B_j^3
    + \ep_n(g,t),
  \]
where
  \[
  \ep_n(g,t) = \ga\sum_{j=1}^{\flr{nt}}g^{(5)}(\be_j)\De B_j^5
    + \sum_{j=1}^{\flr{nt}}h(B(t_{j-1}),B(t_j))\De B_j^6
    + g(B(t)) - g(B_n(t)).
  \]
It will therefore suffice to show that $\ep_n(g,t)\to 0$ ucp.

By the continuity of $g$ and $B$, $g(B(t),t) - g(B_n(t),\flr{nt}/n) \to 0$ uniformly on compacts, with probability one. By Lemma \ref{L:5int_conv}, since $g^{(5)}\in C^1(\RR)$, $\ga\sum_{j=1}^{\flr{nt}} g^{(5)}(\be_j)\De B_j^5\to 0$ ucp. It remains only to show that
  \begin{equation}\label{sextic0ucp}
  \sum_{j=1}^{\flr{nt}} h(B(t_{j-1}),B(t_j))\De B_j^6
    \to 0 \quad\text{ucp}.
  \end{equation}
Fix $T>0$. Let $\{n(k)\}_{k=1}^\infty$ be an arbitrary sequence of positive integers. By Theorem \ref{T:sextic_var}, we may find a subsequence $\{m(k)\}_{k=1}^ \infty$ and a measurable subset $\Om^*\subset\Om$ such that $P(\Om^*) =1$, $t \mapsto B(t,\om)$ is continuous for all $\om\in\Om^*$, and
  \begin{equation}\label{Taylor_expan.1}
  \sum_{j=1}^{\flr{m(k)t}}\De B_{j,m(k)}(\om)^6 \to 15t,
  \end{equation}
as $k\to\infty$ uniformly on $[0,T]$ for all $\om\in\Om^*$. Fix $\om\in\Om^*$. We will show that
  \[
  \sum_{j=1}^{\flr{m(k)t}}
    h(B(t_{j-1}^{m(k)},\om),B(t_j^{m(k)},\om))\De B_{j,m(k)}(\om)^6
    \to 0,
  \]
as $k\to\infty$ uniformly on $[0,T]$, which will complete the proof.

For this, it will suffice to show that
  \[
  \sum_{j=1}^{\flr{m(k)T}}
    |h(B(t_{j-1}^{m(k)},\om),B(t_j^{m(k)},\om))|\De B_{j,m(k)}(\om)^6
    \to 0,
  \]
as $k\to\infty$. We begin by observing that, by \eqref{Taylor_expan.1}, there exists a constant $L$ such that $\sum_{j=1}^{\flr{m(k)T}}\De B_{j,m(k)}(\om)^6 < L$ for all $k$. Now let $\ep>0$. Since $g$ has compact support, $g$ is uniformly continuous. Hence, there exists $\de>0$ such that $|b-a|<\de$ implies $|h(a,b)|<\ep/L$ for all $t$. Moreover, there exists $k_0$ such that $k\ge k_0$ implies $|\De B_{j,m(k)} (\om)|<\de$ for all $1\le j\le\flr{m(k)T}$. Hence, if $k\ge k_0$, then
  \[
  \sum_{j=1}^{\flr{m(k)T}}
    |h(B(t_{j-1}^{m(k)},\om),B(t_j^{m(k)},\om))|\De B_{j,m(k)}(\om)^6
    < \frac{\ep}L\sum_{j=1}^{\flr{m(k)T}}\De B_{j,m(k)}(\om)^6
    < \ep,
  \]
which completes the proof. \qed

\begin{corollary}\label{C:Taylor_expan}
If $g\in C^6(\RR)$ has compact support, then $I_n(g',B,t)\approx X_n(t)$, where for any $T>0$,
  \[
  \sup_n\sup_{t\in[0,T]} E|X_n(t)|^2 < \infty.
  \]
\end{corollary}

\pf This follows immediately from Lemma \ref{L:Taylor_expan} and Theorem \ref{T:3int_mom.N}. \qed

\begin{lemma}\label{L:I_n_rel_cpct}
If $g\in C^6(\RR)$ has compact support, then $\{I_n(g',B)\}$ is relatively compact in $D_\RR[0,\infty)$.
\end{lemma}

\pf Define
  \[
  \begin{split}
  X_n(t) &:= \frac1{12}\sum_{j=1}^{\flr{nt}}
    g'''(\be_j)\De B_j^3,\\
  Y(t) &:= g(B(t)) - g(B(0))\\
  \ep_n(t) &:= I_n(g',B,t) - Y(t) - X_n(t).
  \end{split}
  \]
Since $(x,y,z)\mapsto x+y+z$ is a continuous function from $D_{\RR^3}[0, \infty)$ to $D_{\RR}[0,\infty)$, it will suffice to show that $\{(X_n,Y, \ep_n)\}$ is relatively compact in $D_{\RR^3}[0,\infty)$. By Lemma \ref{L:Taylor_expan}, $\ep_n\to0$ ucp, and therefore in $D_\RR[0,\infty)$. Hence, by Lemma \ref{L:prod_space}, it will suffice to show that $\{X_n\}$ is relatively compact in $D_\RR[0,\infty)$.

For this, we apply Theorem \ref{T:moment_cond} with $\be=4$. Fix $T>0$ and let $0\le s<t\le T$. Let $c=\flr{ns}$ and $d=\flr{nt}$. Note that $q(a+b)^4\le C(|a|^2 +|b|^4)$. Hence, since
$g$ has compact support and, therefore, $g'''$ is bounded,
  \[
  \begin{split}
  E[q(X_n(t) - X_n(s))^4] &= E\bigg[q\bigg(\frac1{12}
    \sum_{j=c+1}^d g'''(\be_j)\De B_j^3\bigg)^4\bigg]\\
  &\le CE\bigg|\sum_{j=c+1}^d
    (g'''(\be_j) - g'''(\be_c))\De B_j^3\bigg|^2
    + CE\bigg|\sum_{j=c+1}^d g'''(\be_c)\De B_j^3\bigg|^4\\
  &\le CE\bigg|\sum_{j=c+1}^d
    (g'''(\be_j) - g'''(\be_c))\De B_j^3\bigg|^2
    + CE\bigg|\sum_{j=c+1}^d \De B_j^3\bigg|^4.
  \end{split}
  \]
Since $g'''\in C^3(\RR)$, we may apply Theorems \ref{T:3int_mom} and \ref{T:4th_mom}, which give
  \[
  E[q(X_n(t) - X_n(s))^4] \le C|t_d - t_c|^{4/3} + C|t_d - t_c|^2
    \le C\bigg(\frac{\flr{nt} - \flr{ns}}n\bigg)^{4/3},
  \]
which verifies condition \eqref{moment_cond} of Theorem \ref{T:moment_cond}. As above,
  \[
  \begin{split}
  E|X_n(T)|^2 &\le CE\bigg|\sum_{j=1}^{\flr{nT}}
    (g'''(\be_j) - g'''(\be_c))\De B_j^3\bigg|^2
    + CE\bigg|\sum_{j=1}^{\flr{nT}} \De B_j^3\bigg|^2\\
  &\le CE\bigg|\sum_{j=1}^{\flr{nT}}
    (g'''(\be_j) - g'''(\be_c))\De B_j^3\bigg|^2
    + C\bigg(E\bigg|
    \sum_{j=1}^{\flr{nT}}\De B_j^3\bigg|^4\bigg)^{1/2}\\
  &\le CT^{4/3} + CT.
  \end{split}
  \]
Hence, $\sup_nE|X_n(T)|^2<\infty$. By Theorem \ref{T:moment_cond}, $\{X_n\}$ is relatively compact, completing the proof. \qed

\begin{lemma}\label{L:Hermite_expan}
If $g\in C^9(\RR)$ has compact support, then
  \[
  I_n(g',B,t) \approx g(B(t)) - g(B(0))
    + \frac1{12\sqrt{n}}\sum_{j=1}^{\flr{nt}}
    \frac{g'''(B(t_{j-1})) + g'''(B(t_j))}2 h_3(n^{1/6}\De B_j).
  \]
\end{lemma}

\pf Using the Taylor expansions in the proof of Lemma \ref{L:Taylor_expan}, together with Lemma \ref{L:5int_conv}, we have
  \[
  \sum_{j=1}^{\flr{nt}}g'''(\be_j)\De B_j^3
    \approx \sum_{j=1}^{\flr{nt}}
    \frac{g'''(B(t_{j-1})) + g'''(B(t_j))}2\De B_j^3.
  \]
By Lemma \ref{L:Taylor_expan}, since $h_3(x)=x^3-3x$, it therefore suffices to show that
  \[
  n^{-1/3}\sum_{j=1}^{\flr{nt}}
    \frac{g'''(B(t_{j-1})) + g'''(B(t_j))}2\De B_j
    = n^{-1/3}I_n(g''',B,t) \approx 0.
  \]
Since $g'''\in C^6(\RR)$, this follows from Lemma \ref{L:I_n_rel_cpct}, Corollary \ref{C:Taylor_expan}, and Lemma \ref{L:in_prob}. \qed

\bigskip

\noindent{\bf Proof of Theorem \ref{T:main}.} We first assume that $g$ (and also $G$) has compact support. By Lemma \ref{L:Hermite_expan} and Theorem \ref{T:main_fdd}, we need only show that $\{(B,V_n(B),I_n(g,B))\}$ is relatively compact in $D_{\RR^3}[0, \infty)$. By Lemma \ref{L:prod_space}, it will suffice to show that $\{I_n(g,B)\}$ is relatively compact in $D_\RR[0,\infty)$. But this follows from Lemma \ref{L:I_n_rel_cpct}, completing the proof when $g$ has compact support.

Now consider general $g$. Let
  \[
  \ts{\Xi_n = (B, V_n(B), I_n(g,B))
    \quad\text{and}\quad
    \Xi=(B, \cub{B}, \int g(B)\,dB).}
  \]
For $T>0$, define $\Xi_n^T(t)=\Xi_n(t)1_{\{t<T\}}$ and $\Xi^T(t) =\Xi(t)1_{\{t<T\}}$. By (3.5.2) in \cite{EK}, if two c\`adl\`ag functions $x$ and $y$ agree on the interval $[0,T)$, then $r(x,y)\le e^{-T}$, where $r$ is the metric on $D_{\RR^d}[0,\infty)$. Hence, by Lemma \ref{L:truncate}, it will suffice to show that $\Xi_n^T\to\Xi^T$ in law, where $T>0$ is fixed.

Let $H:D_{\RR^3}[0,\infty)\to\RR$ be continuous and bounded, with $M =\sup|H(x)|$. Define $X_n=H(\Xi_n^T)$ and $X=H(\Xi^T)$, so that it will suffice to show that $X_n\to X$ in law. For each $k>0$, choose $G_k\in C^6(\RR)$ with compact support such that $G_k=G$ on $[-k,k]$. Let $g_k=G_k'$,
  \[
  \ts{\wt\Xi_{n,k} = (B, V_n(B), I_n(g_k,B)),
    \quad\quad
    \wt\Xi_k = (B, \cub{B}, \int g_k(B)\,dB),}
  \]
$X_{n,k}=H(\wt\Xi_{n,k}^T)$ and $Y_k=H(\wt\Xi_k^T)$. Note that $E|X_n - X_{n,k}| \le \de_k$, where
  \[
  \de_k = 2MP\left({\sup_{0\le t\le T}|B(t)| \ge k}\right).
  \]
Also note that that $\de_k\to0$ as $k\to\infty$. Since $G_k$ has compact support, we have already proven that $X_{n,k}\to Y_k$ in law. Hence, by Lemma \ref{L:truncate}, it will suffice to show that $Y_k\to X$ in law. However, it is an immediate consequence of \eqref{main_int_def} that $\wt\Xi_k^T\to\Xi^T$ ucp, which completes the proof. \qed

\bigskip

\noindent{\bf Proof of Theorem \ref{T:joint_conv}.} As in the proof of Theorem \ref{T:main}, $\{(B,V_n(B),J_n)\}$ is relatively compact. Let $(B,X,Y)$ be any subsequential limit. By Theorem \ref{T:signed_cubic}, $X=\ka W$, where $W$ is a standard Brownian motion, independent of $B$. Hence, $(B,X,Y)=(B,\cub{B},Y)$. Fix $j\in\{1,\ldots,k\}$. By Theorem \ref{T:main}, $(B,\cub{B},Y_j)$ has the same law as $(B,\cub{B},\int g_j(B)\,dB)$. Using the general fact we observed at the beginning of the proof of Theorem \ref{T:main_fdd}, together with \eqref{main_int_def} and the definition of $\cub{B}$, this implies $(\int g_j(B)\,dB,Y_j)$ and $(\int g_j(B)\,dB,\int g_j(B)\, dB)$ have the same law. Hence, $Y_j=\int g_j(B)\,dB$ a.s., so $(B,X,Y) =(B,\cub{B},J)$. \qed



\end{document}